\documentstyle[verbatim,12pt,leqno]{amsart}
 
\let\mathcal\cal
\newtheorem{theorem}{Theorem}[section]

\newtheorem{lemma}[theorem]{Lemma}
\newtheorem{corollary}[theorem]{Corollary}

\newtheorem{claim}{Claim}[lemma]
\newtheorem{problem}{Problem}
\newtheorem{case}{Case}
\newtheorem{subcase}{Subcase}[case]
\newtheorem{Claim}{Claim}

\theoremstyle{definition}

\newtheorem{definition}[theorem]{Definition}

\theoremstyle{remark}

\newcommand{\proof}{\begin{pf}}
\newcommand{\Proof}[1]{\begin{pf*}{Proof of #1}}
\newcommand{\eproof}{\end{pf}}
\newcommand{\Eproof}{\end{pf*}}
\newcommand{\sproof}[1]{\begin{pf*}{#1}}
\newcommand{\esproof}{\end{pf*}}

\newcommand{\arablabel}{
          \renewcommand{\labelenumi}{{\rm (\arabic{enumi})}}
          \renewcommand{\theenumi}{{\rm (\arabic{enumi})}}
          \renewcommand{\labelenumii}{{\rm (\arabic{enumii})}}
          \renewcommand{\theenumii}{{\rm (\arabic{enumii})}}
                    }

\newcommand{\alabel}{
          \renewcommand{\labelenumi}{{\rm (\alph{enumi})}}
          \renewcommand{\theenumi}{{\rm (\alph{enumi})}}
          \renewcommand{\labelenumii}{{\rm (\alph{enumii})}}
          \renewcommand{\theenumii}{{\rm (\alph{enumii})}}
                    }

\newcommand{\rlabel}{
          \renewcommand{\labelenumi}{{\rm (\roman{enumi})}}
          \renewcommand{\theenumi}{{\rm (\roman{enumi})}}
          \renewcommand{\labelenumii}{{\rm (\roman{enumii})}}
          \renewcommand{\theenumii}{{\rm (\roman{enumii})}}
                    }

\def\myheads#1;#2;{
\pagestyle{myheadings}
\markboth{{\sc\hfill #1\hfill\protect\makebox[0cm][r]{\rm\today}}}
{{\sc\protect\makebox[0cm][l]{\rm\today}\hfill #2\hfill}}
}

\newcommand{\bcal}{{\mathcal B}}

\newcommand{\dcal}{{\mathcal D}}

\newcommand{\fcal}{{\mathcal F}}
\newcommand{\gcal}{{\mathcal G}}
\newcommand{\hcal}{{\mathcal H}}

\newcommand{\kcal}{{\mathcal K}}

\newcommand{\pcal}{{\mathcal P}}

\newcommand{\setm}{\setminus}
\newcommand{\empt}{\emptyset}
\newcommand{\subs}{\subset}

\newcommand{\oo}{{{\omega}_1}}
\newcommand{\rest}{\lceil}
\newcommand{\dom}{\operatorname{dom}}

\def\<{\left\langle}
\def\>{\right\rangle}

\def\OO{{\omega}}
\def\oo{\omega_1}
\def\br#1;#2;{\bigl[ {#1} \bigr]^ {#2} }
\def\bc#1;#2;{\bigl( {#1} \bigr)^ {#2} }

\def\ooseq#1;#2;{\< {#1}_{#2}:{#2}<\oo\>}
\def\ooset#1;#2;{\{ {#1}_{#2}:{#2}<\oo\}}
\def\seq#1;#2;#3;{\< {#1}_{#2}:{#2}<#3\>}
\def\set#1;#2;#3;{\{ {#1}_{#2}:{#2}<#3\}}
\def\oseq#1;#2;{\< {#1}_{#2}:{#2}<\OO\>}
\def\oset#1;#2;{\{ {#1}_{#2}:{#2}<\OO\}}
\def\oosequ#1;#2;{\< {#1}^{#2}:{#2}<\oo\>}
\def\oosetu#1;#2;{\{ {#1}^{#2}:{#2}<\oo\}}
\def\sequ#1;#2;#3;{\< {#1}^{#2}:{#2}<#3\>}
\def\setu#1;#2;#3;{\{ {#1}^{#2}:{#2}<#3\}}
\def\osequ#1;#2;{\< {#1}^{#2}:{#2}<\OO\>}
\def\osetu#1;#2;{\{ {#1}^{#2}:{#2}<\OO\}}

\def\force{\raisebox{1.5pt}{\mbox {$\scriptscriptstyle\|$}}
\mbox{$\!\mbox{---}$}}

\newcommand{\fn}{\operatorname{Fn}}
\newcommand{\HH}{\operatorname{H}}

\def\to{\longrightarrow}

\newcommand{\leo}{\operatorname{<_{\text{\rm{On}}}}}

\newcommand{\newcases}{\setcounter{case}{0}}
\newcommand{\newClaim}{\setcounter{Claim}{0}}

\def\fin#1;{\br #1;<{\omega};}

\def\ahii#1;{\<a^{#1},h^{#1},i^{#1}\>}

\def\aii#1;{a^{#1}}
\def\apx{\aii{p_{\xi}};}

\def\ap{a^p}
\def\app{a^{p'}}
\def\hp{{h}^p}
\def\ip{i^p}
\def\aq{a^q}
\def\hq{{h}^q}
\def\iq{i^q}
\def\ar{a^r}
\def\hr{{h}^r}
\def\ir{i^r}
\def\as{a^s}
\def\at{a^t}
\def\hs{{h}^s}
\def\is{i^s}
\def\itt{i^t}
\newcommand{\tght}{\operatorname{t}}

\newcommand{\z}{\operatorname{z}}
\newcommand{\clf}{\operatorname{cl}_f}

\newcommand{\supp}{\operatorname{supp}}
\newcommand{\oot}{{\omega}_2}
\theoremstyle{plain}
  
\newtheorem{Prop}{Proposition}  

\def\adot{\dot{A}}

\def\yydot{\dot{y}}

\def\udot{\dot{U}}
\def\rdot{\dot{R}}
\def\rrdot{\dot{r}}

\def\qqdot{\dot{q}}
\def\rrdot{\dot{r}}
\def\ydot{\dot{Y}}
\def\zzdot{\dot{z}}

\def\hhh#1;#2;#3;{{h}^{#1}[#2,#3]}
\def\hh#1;#2;{{h}^{#1}[#2]}
\def\gg#1;#2;{{g}^{#1}[#2]}
\def\UU#1;#2;{\operatorname{U}(#1,#2)}
\def\uu#1;#2;{{u}(#1,#2)}
\def\uuu#1;#2;#3;{{u}^{#1}(#2,#3)}

\def\hhr#1{{h}^r[#1]}
\def\hhs#1{{h}^s[#1]}
\def\hht#1{{h}^t[#1]}

\def\HHH#1{\operatorname{H}[#1]}
\def\KKK#1{\operatorname{K}(#1)}
\def\HHHH#1;#2;{\operatorname{H}[#1,#2]}

\myheads How to force$\dots$;How to force$\dots$;

\author{I. Juh\'asz}

\address{Mathematical Institute of the Hungarian Academy of Sciences}

\email{juhasz@@math-inst.hu}

\author{L. Soukup}
\address{Mathematical Institute of the Hungarian Academy of Sciences}

\email{soukup@@math-inst.hu}

\title{How to force a countably tight, initially $\oo$-compact and 
non-compact space?}

\thanks{The preparation of this paper was supported by the 
Hungarian National Foundation for Scientific Research grant no. 16391.}



\begin{document}
\maketitle

\begin{center}
\today
\end{center}
\begin{abstract}
Improving a  result  of M. Rabus we force a normal, locally compact, 0-dimensional,
Frechet-Uryson,  
initially $\oo$-compact and  non-compact space $X$ of size 
$\oot$ having the following property: for every open 
(or closed) set $A$ in $X$ we have  $|A|\le\oo$ or  $|X\setm A|\le\oo$.
\end{abstract}

\section{Introduction}\label{sc:int}

E. van Douwen and, independently, A. Dow \cite{Dow}
have observed that under CH an initially $\oo$-compact $T_3$ space
of countable tightness is compact.
(A space $X$ is initially ${\kappa}$-compact if any open cover
of $X$ of size $\le {\kappa}$ has a finite subcover, or
equivalently any subset of $X$ of size $\le {\kappa}$ has
a complete accumulation point). Naturally, the question arose whether
CH is needed here, i.e. whether the same is provable just in ZFC. 
The question became even more intriguing when in \cite{FNy} 
D. Fremlin and P. Nyikos proved the same result from PFA.
Quite recently, A. V. Arhangel'ski\u\i has devoted the paper 
\cite{Arh} to this problem, in which he has raised many
related problems as well.

In \cite{Ra}  M. Rabus has answered 
the question of van Douwen and Dow in the negative.
He constructed by forcing a Boolean algebra B such that the Stone space 
$St(B)$ includes a counterexample $X$ of size $\oot$ to the
van Douwen--Dow question, in fact $St(B)$ is the one point 
compactification of $X$, hence $X$ is also locally compact. 
The forcing used by Rabus is closely related to 
the one due to J. Baumgartner and S. Shelah in \cite{BS}, which
had been used to construct a thin very tall superatomic Boolean algebra.
In particular, Rabus makes use of a so-called $\Delta$-function $f$
(which was also used and introduced in \cite{BS}) with some extra 
properties that are satisfied if $f$ is obtained by the original,
rather sophisticated forcing 
argument of  Shelah from \cite{BS}.

In this paper we give an alternative forcing construction of counterexamples
to the van Douwen--Dow question, which we think is simpler,
more direct and more intuitive than the one in \cite{Ra}.
First of all, we directly force a topology ${\tau}_f$ on $\oot$
that yields an example from a $\Delta$-function (with no extra properties)
in the ground model which also satisfies CH.
There is a wide variety of such ground models since they are easily 
obtained when one forces a $\Delta$-function or because 
 $\raisebox{3pt}{\fbox{}}_{\oo}$  implies the existence of a 
$\Delta$-function (cf. \cite{BS}).

Let us recall  the definition of the $\Delta$-functions 
 from  \cite{BS}.
\begin{definition}\label{df:delta-function}
Let  $f:\br \oot;2;\to \br \oot;\le{\omega};$ be a function with
$f\{{\alpha},{\beta}\}\subs {\alpha}\cap {\beta}$ for 
$\{{\alpha},{\beta}\}\in\br \oot;2;$.
 (1) We say that two finite subsets $x$ and $y$ of $\oot$ are 
{\em good for $f$} provided that for  ${\alpha}\in x\cap y$,
${\beta}\in x\setm y$ and ${\gamma}\in y\setm x$ we always have 
\begin{enumerate}\alabel
\item ${\alpha}<{\beta},{\gamma}$ $\Longrightarrow$ 
${\alpha}\in f\{{\beta},{\gamma}\}$,
\item ${\alpha}<{\beta}$ $\Longrightarrow$ 
$f\{{\alpha},{\gamma}\}\subs f\{{\beta},{\gamma}\}$,
\item ${\alpha}<{\gamma}$ $\Longrightarrow$
$f\{{\alpha},{\beta}\}\subs f\{{\gamma},{\beta}\}$.
\end{enumerate} 
(2) We say that $f$ is a {\em ${\Delta}$-function} 
if every uncountable
family of finite subsets of $\oot$ contains two sets $x$ and $y$
which are good for $f$.
\end{definition}

Both in \cite{BS} and \cite{Ra} the main use of the 
$\Delta$-function $f$ is to suitably restrict the partial order of
finite approximations to a structure on $\oot$ so as to become
c.c.c. This we do as well, but in the proof
of the countable compactness of ${\tau}_f$ we also need the following
simple result that yields an additional property of $\Delta$-functions
provided CH  also holds. In fact, only property \ref{df:delta-function}.(a)
is needed for this.

\begin{lemma}\label{lm:big_int}
Assume that CH holds,  $f$ is a $\Delta$-function, 
 $\{c_{\alpha}:{\alpha}<\oot\}$ are pairwise disjoint finite subsets
of $\oot$ and $B\in\br \oot;{\omega};$. Then for each  $n\in {\omega}$ 
there are distinct ordinals 
${\alpha}_0,{\alpha}_1,\dots,{\alpha}_{n-1}\in\oot$
such that 
$$
B\subs \bigcap\{f({\xi},{\eta}):
{\xi}\in c_{{\alpha}_i},{\eta}\in c_{{\alpha}_j}, i<j<n\}.
$$
\end{lemma}

\proof
We can assume that $\sup B<\min c_{\alpha}$ for each ${\alpha}<\oot$.
Denote by S(n)  the statement of the lemma for $n$. 
We prove S(n) by induction on $n$. 
The first non-trivial case is n=2. 
Assume indirectly that S(2)  fails. 
Then for each ${\alpha}<{\beta}<\oot$ 
there is  $b_{{\alpha},{\beta}}\in B$  such that 
$b_{{\alpha},{\beta}}\notin f({\xi},{\eta})$ for some 
${\xi}\in c_{\alpha}$ and ${\eta}\in c_{\beta}$.
By CH the Erd\H os- Rado partition theorem \cite{Ku} 
has the form $\oot\to(\oo)^2_{\omega}$, thus there are $I\in\br\oot;\oo;$
and $b\in B$ such that for each  ${\alpha}\ne{\beta}\in I$ we have 
$b\notin f({\xi},{\eta})$ for some 
${\xi}\in c_{\alpha}$ and ${\eta}\in c_{\beta}$.

Let $d_{\alpha}=c_{\alpha}\cup\{b\}$ for ${\alpha}\in I$.
Since $f$ is a $\Delta$-function there are ${\alpha}\ne {\beta}\in I$ such
that $d_{\alpha}$ and $d_{\beta}$ are good for $f$.
But $d_{\alpha}\cap d_{\beta}=\{b\}$ and $b<\min c_{\alpha},\min c_{\beta}$,
so by \ref{df:delta-function}.(a) we have 
$b\in f({\xi},{\eta})$ for each ${\xi}\in c_{\alpha}$ and
${\eta}\in c_{\beta}$ contradicting the choice of 
$b=b_{{\alpha},{\beta}}$. Thus S(2) holds.

Assume now that S(n) holds for some $n\ge 2$ and we prove
S(2n). 
Applying 
S(n) $\oot$-many times for $B$ and suitable final segments of
$\{c_{\alpha}:{\alpha}<\oot\}$ we can obtain $\oot$-many pairwise disjoint 
$n$-element sets
$\{{\alpha}^{\nu}_0,{\alpha}^{\nu}_1,\dots,{\alpha}^{\nu}_{n-1}\}\subs \oot$
 such that for each ${\nu}<\oot$
$$
B\subs \bigcap\{f({\xi},{\eta}):
{\xi}\in c_{{\alpha}^{\nu}_i},{\eta}\in c_{{\alpha}^{\nu}_j}, i<j<n\}.
$$
%We can assume that ${\alpha}^{\nu}_i<{\alpha}^{\mu}_j$ whenever ${\nu}<{\mu}$.
Let $d_{\nu}=\bigcup\{c_{{\alpha}^{\nu}_i}:i<n\}$ for ${\nu}<\oot$. 
Applying  S(2) for $B$ and the sequence 
$\{d_{\nu}:{\nu}<\oot\}$ we  get ordinals
${\nu}<{\mu}<\oot$ such that 
$B\subs f({\xi},{\eta})$ for all ${\xi}\in d_{\nu}$ and ${\eta}\in d_{\mu}$.
In other words, if $i,j<n$, ${\xi}\in c_{{\alpha}^{\nu}_i}$ and 
${\eta}\in c_{{\alpha}^{\mu}_j}$ then $B\subs f({\xi},{\eta})$.
Therefore the set 
$\{{\alpha}^{\nu}_0,{\alpha}^{\nu}_1,\dots,{\alpha}^{\nu}_{n-1},
{\alpha}^{\mu}_0,{\alpha}^{\mu}_1,\dots,{\alpha}^{\mu}_{n-1}\}$ 
witnesses  S(2n). 
\eproof

The following, even simpler, result about arbitrary functions
$f:\br \oot;2;\to \br \oot;\le{\omega};$ with
$f\{{\alpha},{\beta}\}\subs {\alpha}\cap {\beta}$ for 
$\{{\alpha},{\beta}\}\in\br \oot;2;$ will also be needed.

\begin{lemma}\label{lm:closure}
If $f$ is a function as above then 
for each $K,K'\in \br\oot;\le{\omega};$ there is a countable set  
$\clf(K,K')\subs \oot$ such that 
\begin{enumerate}\alabel
\item $K\subs \clf(K,K')$, $\sup K=\sup \clf(K,K')$,
%\item $\forall {\xi}\ne {\eta}\in \clf(K,K')$ 
%$f({\xi},{\eta})\subs \clf(K,K')$,
\item $\forall {\xi}\in \clf(K,K')$ $\forall {\eta}\in \clf(K,K')\cup K'$ 
$\left(f\{{\xi},{\eta}\}\subs \clf(K,K')\right)$.
\end{enumerate}
\end{lemma}

\proof
Let $\KKK 0=K$, $\KKK {n+1}=\KKK n\cup\bigcup\{f\{{\xi},{\eta}\}:
{\xi}\in \KKK n, {\eta}\in \KKK n\cup K'\}$ and $\clf(K,K')=\bigcup_{n<{\omega}}\KKK n$. 
\eproof

The topology ${\tau}_f$ that we will construct on $\oot$ is right 
separated (in the natural order of $\oot$) and
is also locally compact and 0-dimensional. Thus for each ${\alpha}\in\oot$
one can fix a compact (hence closed) and open neighbourhood
$\HH({\alpha})$ of ${\alpha}$ such that $\max \HH({\alpha})={\alpha}$. Conversely, if we can fix for each ${\alpha}\in\oot $ 
such a right-separating compact open neighbourhood $\HH({\alpha})$ 
then the family $\{\HH({\alpha}):{\alpha}<\oot\}$
determines the whole topology ${\tau}$ on $\oot$. In fact, using the 
notation
$\UU {\alpha};b;=\HH({\alpha})\setm \bigcup\{\HH({\beta}):{\beta}\in b\}$,
it is easy to check that for each ${\alpha}\in\oot$ the family 
$\bcal_{\alpha}=\{\UU {\alpha};b;:b\in\br {\alpha};<{\omega};\}$
is a ${\tau}$-neighbourhood base of ${\alpha}$. Therefore, our notion of forcing consists of finite approximations
to a family $\hcal=\{\HH({\alpha}):{\alpha}<\oot\}$ like above.

Now, if $\hcal=\{\HH({\alpha}):{\alpha}<\oot\}$ is as required and
${\beta}<{\alpha}<\oot$ then either (i) ${\beta}\in \HH({\alpha})$
or (ii) ${\beta}\notin \HH({\alpha})$. If (i) holds then
$\HH({\beta})\setm \HH({\alpha})$, if (ii) holds then
$\HH({\beta})\cap \HH({\alpha})$ is a compact open subset of ${\beta}$,
hence there is a finite subset of ${\beta}$, call it $i\{{\alpha},{\beta}\}$,
such that this set is covered by 
$\HHH{i\{{\alpha},{\beta}\}}=
\bigcup\{\HH({\gamma}):{\gamma}\in i\{{\alpha},{\beta}\}\}$
It may come as a surprise, but the existence of such a function
$i$ is also sufficient to insure that the collection $\hcal$
be as required. More precisely, we have the following result.

\begin{definition}\label{df:tau_h}
If $\hcal=\{\HH({\alpha}):{\alpha}\in\oot\}$ is  a family of
subsets of $\oot$ such that 
$\max\HH({\alpha})={\alpha}$ for each ${\alpha}\in\oot$ then we denote
by ${\tau}_{\hcal}$ the topology on $\oot$ generated by
$\hcal\cup\{\oot\setm H:H\in\hcal\}$ as a subbase. 
\end{definition}
Clearly,
${\tau}_{\hcal}$ is a 0-dimensional, Hausdorff and right separated
topology in which the elements of $\hcal$ are clopen.
\begin{theorem}\label{tm:loc-comp}
Assume that $\hcal$ is as  in definition
\ref{df:tau_h} above  and there is a function
$i:\br \oot;2;\to \br \oot;<{\omega};$ satisfying 
$i\{{\alpha},{\beta}\}\subs {\alpha}\cap {\beta}$ for each 
$\{{\alpha},{\beta}\}\in \br \oot;2;$ such that if ${\beta}<{\alpha}$
then ${\beta}\in\HH({\alpha})$ implies 
$\HH({\beta})\setm \HH({\alpha})\subs \HHH {i\{{\alpha},{\beta}\}}$ and
${\beta}\notin\HH({\alpha})$ implies 
$\HH({\beta})\cap \HH({\alpha})\subs \HHH {i\{{\alpha},{\beta}\}}$.
Then each $\HH({\alpha})$ is compact in the topology ${\tau}_{\hcal}$,
hence ${\tau}_{\hcal}$ is locally compact. 
\end{theorem}

\proof
We do induction on ${\alpha}\in\oot$. Assume that for each 
${\beta}\in{\alpha}$ we know $\HH({\beta})$ is compact in ${\tau}_{\hcal}$.
By Alexander's subbase lemma 
%(see \cite{Alex}) 
it suffices to show that any 
cover  $\kcal$ of $\HH({\alpha})$ by members of $\hcal$ and their complements
has a finite subcover. Let $K\in\kcal$ be such that 
${\alpha}\in K$. If $K=\HH({\gamma})$ then ${\alpha}\le {\gamma}$.
The case ${\alpha}={\gamma}$ is trivial so assume ${\alpha}<{\gamma}$.
But then $\HH({\alpha})\setm K\subs \HHH{i\{{\alpha},{\gamma}\}}$ and
by our inductive hypothesis $\HH({\beta})$ is compact for each 
${\beta}\in i\{{\alpha},{\gamma}\}$ hence so is 
$\HH({\beta})\cap \HH({\alpha})\setm K$ being closed in $\HH({\beta})$.
Therefore $\HH({\alpha})\setm K$ is compact and so some finite 
$\kcal_0\subs\kcal$ covers it. But then $\kcal_0\cup\{K\}$ covers
$\HH({\alpha})$, hence we are done. A similar argument works if 
$K=\oot\setm \HH({\gamma})$, then using 
$\HH({\alpha})\cap \HH({\gamma})\subs \HHH{i\{{\alpha},{\gamma}\}}$
if ${\alpha}<{\gamma}$ (or the compactness of $\HH({\gamma})$ if 
${\gamma}<{\alpha}$). 
\eproof

It is now very natural to try to force a generic 0-dimensional,
locally compact and right separated topology on $\oot$
by finite approximations (or pieces of information) of $\hcal$ and $i$.
As was already mentioned, the $\Delta$-function $f$ comes into the 
picture when one wants to make this forcing c.c.c. The technical details
of this are done in section \ref{sc:forcing}.

We call the family $\hcal$ {\em coherent} if ${\beta}\in \HH({\alpha})$
implies $\HH({\beta})\subs \HH({\alpha})$. Clearly, this makes 
things easier because then $\HH({\beta})\setm \HH({\alpha})=\empt$,
hence there is no problem covering it, the requirement
on $i$ is only that if ${\beta}\notin \HH({\alpha})$ and ${\beta}<{\alpha}$
then $\HH({\beta})\cap \HH({\alpha})\subs \HHH{i\{{\alpha},{\beta}\}}$.
The original forcing of Baumgartner and Shelah from \cite{BS}
(when translated to scattered, i.e. right separated, locally
compact spaces rather than superatomic Boolean algebras) 
actually produced such a
coherent family $\hcal$. This is interesting because if $\hcal$
is coherent and ${\tau}_{\hcal}$ is separable, which we have almost
automatically if $\hcal$ is obtained generically, then 
${\tau}_{\hcal}$ is also countably tight!

\begin{theorem}\label{tm:coherent}
If there is a coherent family $\hcal$ of right separating compact open sets
for a separable topology ${\tau}$ on $\oot$ then 
$\tght(\<\oot,{\tau}\>)={\omega}$.
\end{theorem}

\proof
Let $X=\<\oot,{\tau}\>$. Then for each ${\alpha}\in\oot$ we have 
$\tght({\alpha},X)=\tght({\alpha},\HH({\alpha}))$,
hence it suffices to prove $\tght({\alpha},\HH({\alpha}))={\omega}$.
We do this by induction on ${\alpha}$. So assume it for each 
${\beta}<{\alpha}$.
If we had $\tght({\alpha},\HH({\alpha}))=\oo$ then $\HH({\alpha})$
would contain a free sequence $S=\{x_{\nu}:{\nu}<\oo\}$ of length
$\oo$. $S$ must converge to ${\alpha}$ since for each 
${\beta}<{\alpha}$ we have, by the inductive hypothesis,
$\tght(\HH({\beta}))\le {\omega}$, hence $|\HH({\beta})\cap S|={\omega}$.
Let $F_{\nu}=\overline{\{x_{\mu}:{\mu}<{\nu}\}}$ for each ${\nu}<\oo$.
Since $S$ is free we have ${\alpha}\notin F_{\nu}$ for all
${\nu}<\oo$ and the sequence $\<F_{\nu}:{\nu}<\oo\>$ is (strictly)
increasing. For each ${\nu}<\oo$ there is a finite subset
$b_{\nu}$ of $F_{\nu}$ such that $F_{\nu}\subs \HHH{b_{\nu}}$.
Now, if ${\mu}<{\nu}$ then $b_{\mu}\subs F_{\mu}\subs F_{\nu}\subs
\HHH{b_{\nu}}$ implies that for each ${\beta}\in b_{\mu}$ 
we have $\HH({\beta})\subs \HHH{b_{\nu}}$ by the coherence of $\hcal$,
hence $\HHH{b_{\mu}}\subs \HHH{b_{\nu}}$.
But we have seen above that $|\HHH{b_{\nu}}\cap S|\le {\omega}$ for each 
${\nu}\in\oo$, while of course $S\subs\bigcup\{\HHH{b_{\nu}}:{\nu}<\oo\}$,
hence the $\HHH{b_{\nu}}${\bf 's} yield a strictly increasing 
$\oo$-sequence of clopen sets in a separable space, which is a 
contradiction
completing the proof.
\eproof

Ironically, this general result that gives countable tightness so easily
cannot be used in our construction because we had to abandon the coherency of
$\hcal$ in our effort to insure countable compactness
(implied by the initial $\oo$-compactness) of ${\tau}_{\hcal}$.

We mentioned above that our examples, by genericity, are separable.
But this is not a coincidence. 
It is well-known and very easy to prove that if $X$ is an 
initially $\oo$-compact space then  $\tght(X)\le{\omega}$ implies that 
$X$ has no
uncountable free sequence. (Moreover, if $X$ is $T_3$ the converse 
of this is also true.) 
Hence the following easy, but perhaps not widely known, result immediately
implies that any non-compact, initially $\oo$-compact space of countable
tightness contains a countable subset whose closure is not compact. Thus if
 there is a counterexample to the  van Douwen--Dow question
then there is also a separable one.

\begin{lemma}\label{lm:noncompact_free_sequence}
If $Y$ is a non-compact topological space, then  for
some ordinal $\mu$ the space $Y$ contains a  free
sequence $\{y_\xi:\xi<\mu\}\subs Y$ with non-compact closure.  
\end{lemma}

\proof
If $Y$ is non-compact, then $Y$ has an strictly increasing 
open cover $\{U_\alpha:\alpha<\kappa\}$ for some regular cardinal
$\kappa$. We pick  points $y_\xi\in X$  and ordinals 
$\alpha_\xi<\kappa$ by recursion as follows.
If the closure of the set $Y_\xi=\{y_\eta:\eta<\xi\}$ in $Y$ is compact, 
then pick $\alpha_\xi\in \kappa$ such that 
$\overline{Y_\xi}\subs U_{\alpha_{\xi}}$
and let $y_\xi\in Y\setm U_{\alpha_\xi}$.  

The sequence $\alpha_\xi$ is strictly increasing 
because $y_\eta\in U_{{\alpha}_\xi}\setm U_{{\alpha}_\eta}$ for $\eta<\xi$.
 So for some $\xi\le\kappa$ the closure of $Y_\xi$ is
non-compact. But $Y_\xi$ is also free because for each ${\eta}<{\xi}$ we have
$\overline{Y_\eta}\subs U_{\alpha_\eta}$ and
$(Y_\xi\setm Y_\eta)\cap U_{\alpha_\eta}=\empt$. 
So we are done.
\eproof

Note that under CH the weight of a separable $T_3$ space is 
$\le\oo$, and an initially $\oo$-compact space of weight $\le\oo$
is compact, hence the CH result of van Douwen and Dow is a trivial 
consequence of \ref{lm:noncompact_free_sequence}.   
Arhangel'ski\u\i raised the question, \cite[problem 3]{Arh},
whether in this CH can be weakened to $2^{\omega}<2^{\oo}$? 
We shall answer this question in the negative: 
 theorem \ref{tm:FU} implies
that the existence of a counterexample to the van Douwen--Dow 
question is consistent
with practically any cardinal arithmetic that   violates CH.

In \cite[problem 17]{Arh}  Arhangel'ski\u\i asked if 
it is provable in ZFC that an initially
$\oo$-compact subspace of a $T_3$ space of countable tightness is 
always closed.
(Clearly this is so under CH or PFA, or in general if the answer to 
the van Douwen--Dow question is ``yes''.) In view of our next result both
Rabus' and our spaces give a negative answer to this question. 
More generally we have the following result.
\begin{theorem}
If $X$ is a locally compact counterexample to van Douwen--Dow 
then the one-point
compactification ${\alpha}X=X\cup\{p\}$ of $X$ also has countable tightness.
On the other hand, $X$ is an initially $\oo$-compact non-closed subset of
${\alpha}X$.
\end{theorem}

\proof
Let $A\subs X$ be such that $p\in\overline{A}$ ( i.e. 
$\overline{A}^X$ is not compact). By lemma 
\ref{lm:noncompact_free_sequence} and our preceding  remark then there is a countable 
set $S\subs \overline{A}^X$ such that 
$\overline{S}^X$ is not compact. But by 
$\tght(X)={\omega}$ then there is a countable $T\subs A$ for which
$S\subs\overline{T}^X$, hence $\overline{T}^X$ is non-compact as well,
so $p\in\overline{T}$. Consequently we have $\tght(p,{\alpha}X)={\omega}$
and so $\tght({\alpha}X)={\omega}$.
\eproof

\section{The forcing construction }\label{sc:forcing}

The following notation will be used in the definition of the poset  $P_f$.
Given a function  $h$ and $a\subs \dom(h)$ we write 
$\hh {};{a};=\cup\{h({\xi}):{\xi}\in a\}$.
Given non-empty sets $x$ and $y$ of ordinals with $\sup x\ne\sup y$ let 
\begin{displaymath}
x*y= \left\{ 
\begin{array}{ll}
x\cap y&\mbox{if $\sup x\notin y$ and $\sup y \notin x$,}\\
x\setm y&\mbox{if $\sup x \in y$,}\\
y\setm x&\mbox{if $\sup y \in x$.}
\end{array}
\right.
\end{displaymath}

\begin{definition}\label{df:poset}
For each function $f:\br \oot;2;\to \br\oot;\le{\omega};$
satisfying
$f({\alpha},{\beta})\subs {\alpha}\cap {\beta}$ 
for any  $\{{\alpha},{\beta}\}\in\br \oot;2;$ we
define a poset $P_f=\<P_f,\le\>$ as follows.
The underlying set of  $P_f$ is the family of triples $p=\ahii ;$ for which  
\begin{enumerate}\rlabel
\item $a\in\fin \oot;$, $h:a\to\pcal(a)$ and 
$i:\br a;2;\to\pcal(a)$ are functions, 
\item $\max h({\xi})={\xi}$ for each ${\xi}\in a$ ,
\item  $i({\xi},{\eta})\subs f({\xi},{\eta})$  
for each $\{{\xi},{\eta}\}\in \br a;2;$,
\item $h({\xi})*h({\eta})\subs \hh; {i({\xi},{\eta})};$ 
for each $\{{\xi},{\eta}\}\in \br a;2;$.
\end{enumerate}
We will often write $p=\ahii p;$ for $p\in P_f$.\newline
For $p,q\in P_f$ let $p\le q$ if
$\ap\supset \aq$, $\hp({\xi})\cap \aq= \hq({\xi})$ for
${\xi}\in \aq$, and $\ip\supset \iq$.
If $p\in P_f$, ${\alpha}\in \ap$, $b\subs \ap\cap {\alpha}$,
let us write $\uuu p;{\alpha};b;=\hp({\alpha})\setm\hh p;b;$.
\end{definition}

\begin{lemma}\label{lm:dense}
For each ${\alpha}<\oot$ the set $D_{\alpha}=\{p\in P_f:{\alpha}\in \ap\}$ 
is dense in $P_f$.
\end{lemma}

\proof
Let $q\in P_f$ with ${\alpha}\notin\aq$. Define 
the condition $p\le q$ by the stipulations
$\ap=\aq\cup\{{\alpha}\}$, $\hp({\alpha})=\{{\alpha}\}$, 
$\hp({\xi})=\hq({\xi})$ and $\ip({\alpha},{\xi})=\empt$ for ${\xi}\in\aq$.
Then clearly $p\in D_{\alpha}$. 
\eproof

\begin{definition}
\label{df:xf}
If $\gcal$ is a  $P_f$-generic filter over $V$, in $V[\gcal]$
we can define the topological space
 $X_f[\gcal]=X_f=\<\oot,{\tau}_f\>$ as follows. 
For ${\alpha}\in \oot$ put
$\HH({\alpha})=\bigcup\{\hp({\alpha}):p\in\gcal\land{\alpha}\in\ap\}$,
let $\hcal=\{\HH({\alpha}):{\alpha}<\oot\}$
and  let ${\tau}_f={\tau}_{\hcal}$ as defined in \ref{df:tau_h}, that is, 
${\tau}_f$ is  the topology on $\oot$ generated by
$\hcal\cup\{\oot\setm H:H\in\hcal\}$ as a subbase.
\end{definition}

If $\gcal$ is a $P_f$-generic filter over $V$ then
by lemma \ref{lm:dense} we have
 $\bigcup\{\ap:p\in\gcal\}=\oot$, and for each ${\alpha}<\oot$ 
$\max \HH({\alpha})={\alpha}$ and $\HH({\alpha})$ is clopen in $X_f$. 
Thus  $X_f$ is  0-dimensional and  right separated.
Of course,  neither $f$ nor $i$ is needed for this. As was explained 
in section \ref{sc:int}, we need $f$ to be a $\Delta$-function
in order to make $P_f$ c.c.c 
(which insures that no cardinal is collapsed), and
the function $i$ is used to make $X_f$ also locally compact.

\begin{theorem}\label{tm:main}
If CH holds and $f$ is a $\Delta$-function, then $P_f$ satisfies the c.c.c and
$V^{P_f}\models$ ``$X_f=\<{\omega}_2,{\tau}_f\>$ 
is a 0-dimensional, right separated, locally compact space 
having the following properties:
\begin{enumerate}
\item[(i)] $t(X_f)={\omega}$,
\item[(ii)] $\forall A\in [{\omega}_2]^{{\omega}_1}$ 
$\exists {\alpha}\in {\omega}_2$ $|A\cap H({\alpha})|={\omega}_1$
\item[(iii)] $\forall A\in [{\omega}_2]^{{\omega}}$
\ $(\overline{A}$ is compact 
or  $|{\omega}_2\setminus \overline{A}|\le {\omega}_1)$.'' 
\end{enumerate}
Consequently, in $V^{P_f}$, $X_f$ is a locally compact, normal,
countably tight, initially $\oo$-compact but non-compact space.
\end{theorem}

\Proof{theorem \ref{tm:main}}
To show that $P_f$ satisfies c.c.c we will proceed in the following way.
We first formulate when two conditions $p$ and $p'$ from $P_f$ are called
{\em good twins} (definition \ref{df:good-twins}), then we 
construct the {\em  amalgamation} $r=p+p'$ of $p$ and $p'$ 
(definition \ref{df:amalgamation}) and show that $r$
is a common extension of $p$ and $p'$ in $P_f$. Finally we prove 
in lemma \ref{lm:c.c.c} that every uncountable family of conditions 
contains a couple of  elements which are good twins.

\begin{definition}\label{df:good-twins}
Let $p=\ahii ;$ and $p'=\ahii \prime;$ be from  $P_f$. 
We say that $p$ and $p'$ are 
{\em good twins} provided 
\begin{enumerate}\arablabel
\item $p$ and $p'$ are {\em twins}, i.e., $|a|=|a'|$ and the natural 
order-preserving bijection $e=e_{p,p'}$ between $a$ and
$a'$ is an isomorphism between $p$ and $p'$:
\begin{enumerate}\rlabel
\item  $h'(e({\xi}))=e''h({\xi})$ for each ${\xi}\in a$, 
\item $i'(e({\xi}),e({\eta}))=e''i({\xi},{\eta})$ for each 
$\{{\xi},{\eta}\}\in\br a;2; $ ,
\item  $e({\xi})={\xi}$ for each ${\xi}\in a\cap a'$,
 \end{enumerate} 
\item  \label{delta} $i({\xi},{\eta})=i'({\xi},{\eta})$ 
for each $\{{\xi},{\eta}\}\in \br a\cap a';2;$,
\item $a$ and $a'$ are good for $f$.
\end{enumerate}
\end{definition}
 
Let us remark that, in view of (ii) and (iii), condition \ref{delta} 
can be replaced by 
``$i({\xi},{\eta})\subs a\cap a'$ for each 
$\{{\xi},{\eta}\}\in \br a \cap a';2;$''.

%The next definition play crucial role in the this proof.

\begin{definition}\label{df:amalgamation}
If $p=\ahii ;$ and $p'=\ahii \prime;$ are good twins we  define the {\em  amalgamation 
$r=\<b,g,j\>$ of $p$ and $p'$} as follows:

Let ${b}={a}\cup{a'}$.
For ${\xi}\in \hh {};{a\cap a'};\cup \hh \prime;a\cap a';$ define 
\begin{displaymath}
{\delta}_{\xi}= 
\min\{{\delta}\in {a}\cap{a'}: 
{\xi}\in{h}({\delta})\cup{h'}({\delta}) \}.
\end{displaymath}

Now, for any ${\xi}\in b$ let 
\begin{equation}\tag{$\bullet$}\label{bul}
{g}({\xi})= \left\{ 
\begin{array}{ll}
{h}({\xi})\cup h'({\xi})
&\mbox{if ${\xi}\in{a}\cap a'$,}\\
{h}({\xi})\cup\{{\eta}\in{a'}\setm{a}:{\delta}_{\eta}\in{h}({\xi})\}
&\mbox{if ${\xi}\in{a}\setm a'$,}\\
{h'}({\xi})\cup\{{\eta}\in{a}\setm{a'}:{\delta}_{\eta}\in{h'}({\xi})\}
&\mbox{if ${\xi}\in{a'}\setm a$.}
\end{array}
\right.
\end{equation}
%(We will see later that ${g}({\xi})$ is well-defined even
%for ${\xi}\in {a}\cap{a'}$)
%
Finally for $\{{\xi},{\eta}\}\in\br b;2;$ let 
\begin{equation}\tag{$\bullet\bullet$}\label{bull}
{j}({\xi},{\eta})= \left\{ 
\begin{array}{ll}
{i}({\xi},{\eta})&\mbox{if ${\xi},{\eta}\in{a}$,}\\
{i'}({\xi},{\eta})&\mbox{if ${\xi},{\eta}\in{a'}$,}\\
f({\xi},{\eta})\cap {b}&\mbox{otherwise.}
\end{array}
\right.
\end{equation}
(Observe that  $j$ is well-defined  because   \ref{df:good-twins}.(2) holds.)
\newline
We will write $r=p+p'$ for the amalgamation
of $p$ and $p'$.
\end{definition}

\begin{lemma}\label{lm:r=p+q}
If $p$ and $p'$ are good twins 
then  their  amalgamation, $r=p+p'$,   is a common extension
of $p$ and $p'$  in $P_f$.
\end{lemma}

\proof
First  we prove two claims.

\begin{claim}\label{lm:1}
Let ${\eta}\in{a}$ and ${\delta}\in{a}\cap{a'}$. Then 
${\eta}\in{h}({\delta})$ if and only if ${\delta}_{\eta}$ is defined  and 
${\delta}_{\eta}\in{h}({\delta})$.
(Clearly, we also have a symmetric version of this statement for
${\eta}\in a'$.)
\end{claim}

\Proof{claim \ref{lm:1}}
Assume first ${\eta}\in {h}({\delta})$. 
Then ${\delta}_{\eta}$ is defined and 
clearly ${\delta}_{\eta}\in h({\delta})$ if ${\delta}_{\eta}={\delta}$. 
So assume ${\delta}_{\eta}\ne {\delta}$.
Since ${i}({\delta}_{\eta},{\delta})\subs {a}\cap{a'}$ and 
$\max {i}({\delta}_{\eta},{\delta})< {\delta}_{\eta}$ we have
${\eta}\notin \hh ;{{i}({\delta}_{\eta},{\delta})};$ by the choice of
${\delta}_{\eta}$. Thus from $p\in P_f$ we have
\begin{equation}
\tag{\dag} {\eta}\notin {h}({\delta}_{\eta})*{h}({\delta}).
\end{equation}

Then  
${h}({\delta}_{\eta})*{h}({\delta})\ne {h}({\delta}_{\eta})\cap {h}({\delta})$
by (\dag). Since ${\eta}\in{h}({\delta})$ implies ${\delta}_{\eta}<{\delta}$,
we actually have 
${h}({\delta}_{\eta})*{h}({\delta})= {h}({\delta}_{\eta})\setm {h}({\delta})$.
Thus  ${\delta}_{\eta}\in{h}({\delta})$ by the definition of the operation
$*$.

On the other hand, if ${\delta}_{\eta}\in{h}({\delta})$, then 
either ${\delta}_{\eta}={\delta}$ or
${h}({\delta}_{\eta})*{h}({\delta})= {h}({\delta}_{\eta})\setm {h}({\delta})$.
Thus ${\eta}\in {h}({\delta})$ because in the latter case
again ${\eta}\notin \hh ;i({\delta}_{\eta},{\delta});$, hence (\dag) holds. 
\Eproof

\begin{claim}
\label{lm:good_def}
If ${\xi}\in a\cap a'$ then 
$$
{g}({\xi})= 
{h}({\xi})\cup\{{\eta}\in{a'}\setm{a}:{\delta}_{\eta}\in{h}({\xi})\}
=
{h'}({\xi})\cup\{{\eta}\in{a}\setm{a'}:{\delta}_{\eta}\in{h'}({\xi})\}.
$$
\end{claim}

\Proof{claim \ref{lm:good_def}}

Conditions \ref{df:good-twins}.1(i) and (iii) imply
$h({\xi})\cap a\cap a'=h'({\xi})\cap a\cap a'$ and so
$$
g({\xi})=h({\xi})\cup h'({\xi})=h({\xi})\cup ((a'\setm a)\cap h'({\xi})).
$$
By claim \ref{lm:1} we have
$$
(a'\setm a)\cap h'({\xi})=
\{{\eta}\in{a'}\setm{a}:{\delta}_{\eta}\in{h'}({\xi})\}.
$$
But by \ref{df:good-twins}.(1) we have
${\delta}_{\eta}\in{h}({\xi})$ iff
${\delta}_{\eta}\in{h'}({\xi})$, hence
 it follows that
$$
g({\xi})=
{h}({\xi})\cup\{{\eta}\in{a'}\setm{a}:{\delta}_{\eta}\in{h}({\xi})\}.
$$
The second equality follows analogously.
\Eproof

Next we check $r\in P_f$.
Conditions  \ref{df:poset}.(i)--(iii) for $r$  are clear  by the construction.
So we should verify  \ref{df:poset}.(iv). 

Let ${\xi}\ne{\eta}\in{b}$ and
${\alpha}\in {g}({\xi})*{g}({\eta})$.
We need to show that ${\alpha}\in \gg ; {{j}({\xi},{\eta})};$.
We will distinguish several cases.
\begin{case}
${\xi},{\eta}\in{a}$ $($ or ${\xi},{\eta}\in{a'}$ $)$.
\end{case} 

Since ${g}({\xi})\cap a= {h}({\xi})$ and ${g}({\eta})\cap a={h}({\eta})$ 
we have $({g}({\xi})*{g}({\eta}))\cap {a}={h}({\xi})* {h}({\eta})$
by the definition of operation $*$. Thus
$({g}({\xi})*{g}({\eta}))\cap a\subs \hh ;{i}({\xi},{\eta});
=\gg ;{j({\xi},{\eta})});\cap a\subs \gg ; {{j}({\xi},{\eta})};$.
So we can assume that ${\alpha}\in{a'}\setm{a}$. 
We know that  ${\delta}_{\alpha}$ is defined 
because ${\alpha}\in g({\xi})\cup g({\eta})$ is also satisfied. Since 
${\alpha}\in g({\xi})$ iff ${\delta}_{\alpha}\in h({\xi})$ 
and ${\alpha}\in g({\eta})$ iff ${\delta}_{\alpha}\in h({\eta})$
by $(\bullet)$ it follows that 
${\delta}_{\alpha}\in{h}({\xi})*{h}({\eta})$. 
Thus there is 
${\nu}\in {i}({\xi},{\eta})$ such that ${\delta}_{\alpha}\in{h}({\nu})$.
But  $i({\xi},{\eta})=j({\xi},{\eta})$ and ${\alpha}\in {g}({\nu})$ 
by $(\bullet)$. 
Thus ${\alpha}\in \gg ; {{j}({\xi},{\eta})};$.

\begin{case}
${\xi}\in{a}\setm {a'}$ and ${\eta}\in{a'}\setm{a}$.
\end{case}

We can assume that ${\alpha}\in{a}$,
since the ${\alpha}\in a'$ case is done symmetrically.

\begin{subcase}
${g}({\xi})*{g}({\eta})={g}({\eta})\setm {g}({\xi})$.
\end{subcase}

Then   ${\alpha}\in\ {g}({\eta})$ and ${\eta}\in\ {g}({\xi})$ so 
${\delta}_{\alpha}$ and ${\delta}_{\eta}$ are both defined and
${\delta}_{\alpha}\in h'({\eta})$, ${\delta}_{\eta}\in h({\xi})$ hold, hence
 ${\alpha}\le{\delta}_{\alpha}<{\eta}<{\delta}_{\eta}<{\xi}$. 
But ${a}$ and ${a'}$ are good for $f$, 
so by \ref{df:delta-function}(a) we have
${\delta}_{\alpha}\in f({\eta},{\xi})\cap {b}={j}({\xi},{\eta})$. Thus 
${\alpha}\in {h}({\delta}_{\alpha})
\subs {g}({\delta}_{\alpha})\subs \gg ; {{j}({\xi},{\eta})};$
which was to be proved.

\begin{subcase}
${g}({\xi})*{g}({\eta})={g}({\xi})\cap {g}({\eta})$ or 
${g}({\xi})*{g}({\eta})={g}({\xi})\setm {g}({\eta})$.
\end{subcase}

 Since now ${\alpha}\in {g}({\xi})*{g}({\eta})\subs {g}({\xi})$, 
by the definition of the  operation $*$ we have 
\begin{equation}
\label{**} |\{{\alpha},{\xi}\}\cap {g}({\eta})|=1.
\end{equation} 
Thus, by the definition of $g({\eta})$, ${\delta}^*=\min\{{\delta}\in{a}\cap{a'}:{\alpha}\in{h}({\delta})\lor
{\xi}\in{h}({\delta})\}$ is well-defined and ${\delta}^*<{\eta}$. 
If ${\delta}^*<{\xi}$ then ${\alpha}\in{h}({\delta}^*)$ and 
by \ref{df:delta-function}(a) we have
${\delta}^*\in f({\xi},{\eta})\cap{b}={j}({\xi},{\eta})$
for ${a}$ and ${a'}$ are good for $f$, 
and so ${\alpha}\in \gg ; {{j}({\xi},{\eta})};$.

Thus we can assume ${\xi}<{\delta}^*$.  
We know that ${\delta}^*={\delta}_{\alpha}$ or ${\delta}^*={\delta}_{\xi}$
by the choice of ${\delta}^*$, but ${\delta}_{\alpha}= {\delta}_{\xi}$
is impossible by (\ref{**}). Thus 
\begin{equation}\label{***}|\{{\alpha},{\xi}\}\cap h({\delta}^*)|=1.
\end{equation}
Since ${\alpha}\in g({\xi})$ implies ${\alpha}\in{h}({\xi})$ 
and we have ${\xi}<{\delta}^*$, 
(\ref{***}) implies ${\alpha}\in {h}({\xi})*{h}({\delta}^*)$ and so
${\alpha}\in \hh ;i({\xi},{\delta}^*);$.
But $i({\xi},{\delta}^*)\subs f({\xi},{\delta}^*)\subs f({\xi},{\eta})$
because $a$ and $a'$ are good for $f$, so \ref{df:delta-function}(b)
or (c) may be applied.
Consequently, we have   $i({\xi},{\delta}^*)\subs j({\xi},{\eta})$ by (\ref{bull}). 
Hence ${\alpha}\in \gg ; {{j}({\xi},{\eta})};$ which was to be proved.

Since we investigated all the cases it follows that 
$r$ satisfies \ref{df:poset}.(iv), that is, $r\in P_f$. 
Since $r\le p,q$ are clear from the construction, the lemma is proved.
\eproof

\begin{lemma}\label{lm:c.c.c}
Every uncountable family $\fcal$ of conditions in $P_f$
contains a couple of  elements which are good twins.
Consequently, $P_f$ satisfies c.c.c.
\end{lemma}

\proof
By standard counting arguments  
$\fcal$ contains an uncountable subfamily $\fcal'$ 
such that every pair $p\ne {p'}\in\fcal'$
satisfies \ref{df:good-twins}.(1)-(2). But $f$ is a $\Delta$-function, so 
there are $p\ne {p'}\in\fcal'$ 
such that $\ap$ and $\app$ are good for $f$, i.e.
$p$ and ${p'}$ satisfies \ref{df:good-twins}.(3), too. In other words, 
$p$ and ${p'}$ are good twins and so  $r=p+{p'}$ is a common extension
 of $p$ and ${p'}$ in $P_f$.
\eproof

Let $\gcal$ be a $P_f$-generic filter over $V$. 
As in definition \ref{df:xf}, let 
$\HH({\alpha})=\bigcup\{\hp({\alpha}):p\in\gcal\land{\alpha}\in\ap\}$
for ${\alpha}\in \oot$,  and let ${\tau}_f$ be the topology
on $\oot$ generated by
$\{\HH({\alpha}):{\alpha}\in\oot\}\cup\{\oot\setm 
\HH({\alpha}):{\alpha}\in\oot\}$ as a subbase.
 Put $i=\bigcup\{\ip:p\in\gcal\}$.

Since  $X_f$ is generated by a clopen subbase and 
$\max(H({\alpha}))={\alpha}$ for each ${\alpha}\in\oot$ by \ref{df:poset}(ii), 
it follows that $X_f$ is 0-dimensional and right separated in its natural
well-order.

The following proposition is clear by \ref{df:poset}.(iv) and by the 
definition of $\HH{}$ and $i$.
\begin{Prop}
$\HH({\alpha})*\HH({\beta})\subs \HHH {i({\alpha},{\beta})}$ for 
${\alpha}<{\beta}<\oot$. So by \ref{tm:loc-comp} 
every $\HH({\alpha})$ is a compact open set in $X_f$.
\end{Prop}

\begin{definition}\label{df:base}
For ${\alpha}\in\oot$ and $b\in\br {\alpha};<{\omega};$ let
$$
\UU{\alpha};b;=\HH({\alpha})\setm \HHH {b}
$$
and let 
$$\bcal_{\alpha}=\{\UU{\alpha};b;:b\in\fin {\alpha};\}.$$
\end{definition}

By theorem \ref{tm:loc-comp}  every $\HH({\alpha})$ is compact and
$\bcal_{\alpha}$ is a neighborhood base of ${\alpha}$ in $X_f$.
Thus $X_f$ is locally compact and the neighbourhood base $\bcal_{\alpha}$
of ${\alpha}$ consists of compact open sets.

Unfortunately, the family $\hcal=\{\HH({\alpha}):{\alpha}<\oot\}$
is not coherent, so we can't apply theorem \ref{tm:coherent}
to prove that $X_f$ is countably tight.
It will however follow from the following result.

\begin{lemma}\label{lm:converge}
In $V^{P_f}$, if a sequence $\{z_{\zeta}:{\zeta}<\oo\}\subs \HH({\beta})$
 converges to ${\beta}$, then there is some ${\xi}<\oo$ such that
${\beta}\in\overline{\{z_{\zeta}:{\zeta}<{\xi}\}}$.
\end{lemma}

\proof
Assume on the contrary that for each ${\xi}<\oo$ we can find a finite
subset $b_{\xi}\subs {\beta}$ such that 
$\{z_{\zeta}:{\zeta}<{\xi}\}\cap \UU {\beta};{b_{\xi}};=\empt$, that is,
$\{z_{\zeta}:{\zeta}<{\xi}\}\subs \HHH {b_{\xi}}$.

Fix now a condition $p\in P_f$ which forces the above described 
situation and decides the value of ${\beta}$.
Then, for each ${\xi}<\oo$ we can choose  a condition $p_{\xi}\le p$
which decides the value of $z_{\xi}$ and $b_{\xi}$. 
We can assume that $\{\apx:{\xi}<\oo\}$ forms a $\Delta$-system
with kernel $D$, $z_{\xi}\in \apx\setm D$ and that 
$z_{\xi}<z_{\eta}$ for ${\xi}<{\eta}<\oo$.

\begin{Claim}
Assume that ${\xi}<{\eta}<\oo$, $p_{\xi}$ and $p_{\eta}$ are good twins and
 $r=p_{\xi}+p_{\eta}$. Then $r\force$``$z_{\xi}\in \HHH {D\cap {\beta}}$''.
\end{Claim}

\Proof{the claim}
Indeed, $z_{\xi}\in \hhr {b_{\eta}}$ because  $r\le p_{\xi}, p_{\eta}$. 
Since $z_{\xi}\in a^{p_{\xi}}\setm a^{p_{\eta}}$,  ($\bullet$) and
claim \ref{lm:good_def}  imply that $z_{\xi}\in \hhr {b_{\eta}}$ holds
if and only if
${\delta}_{z_{\xi}}\in D=a^{p_{\xi}}\cap a^{p_{\eta}}$ is defined and 
${\delta}_{z_{\xi}}\in \hh {p_{\eta}};{b_{\eta}};$. Since
$b_{\eta}\subs {\beta}$, we also have ${\delta}_{z_{\xi}}<{\beta}$
and so $z_{\xi}\in \hh r;{D\cap {\beta}};$. 
\Eproof
Applying lemma \ref{lm:c.c.c} to appropriate final segments of
 $\{p_{\xi}:{\xi}<\oo\}$  we can choose, by induction on ${\mu}<\oo$, 
pairwise different ordinals 
 ${\mu}<{\xi}_{\mu}<{\eta}_{\mu}<\oo$ with ${\eta}_{\mu}<{\xi}_{\nu}$
if ${\mu}<{\nu}$ 
such that $p_{{\xi}_{\mu}}$ and $p_{{\eta}_{\mu}}$ are good twins.
Let $r_{\mu}=p_{{\xi}_{\mu}}+p_{{\eta}_{\mu}}$. Since $P_f$ satisfies c.c.c
there is a condition $q\le p$ such that 
$q\force$`` $|\{{\mu}\in\oo: r_{\mu}\in\gcal\}|=\oo$''.
Thus, by the claim 
$q\force $``$|\{z_{\xi}:{\xi}<\oo\}\cap \HHH {D\cap {\beta}}|=\oo$'',
i.e. the neighbourhood $\UU {\beta};D\cap {\beta};$ of ${\beta}$  
misses uncountably many of the points $z_{\xi}$
which contradicts that ${\beta}$ is the limit  of this sequence.
\eproof

\begin{corollary}\label{lm:to}
$\tght(X_f)={\omega}$.
\end{corollary}

\proof
Assume on the contrary that ${\alpha}\in \oot$ and 
$\tght({\alpha},X_f)=\tght({\alpha},\HH({\alpha}))=\oo$. 
Then there is an ${\omega}$-closed set 
$Y\subs \HH({\alpha})$  that is not closed. 
Since the subspace $H({\alpha})$ is compact and right separated and  
so it is  pseudo-radial
%(see \cite[???]{pseudo})
,  for some regular 
cardinal ${\kappa}$ there is a sequence
 $\{z_{\xi}:{\xi}<{\kappa}\}\subs Y$ which converges 
to some point ${\beta}\in H({\alpha})\setm Y$. Since $Y$ is ${\omega}$-closed
and $|Y|\le |\HH({\alpha})|=\oo$ we have ${\kappa}=\oo$.
By lemma \ref{lm:converge} there is some ${\xi}<\oo$ with
${\beta}\in \overline {\{z_{\zeta}:{\zeta}<{\xi}\}}\subs Y$ contradicting
${\beta}\notin Y$.
\eproof

\begin{lemma}\label{lm:oo-compact}
In $V^{P_f}$, for each uncountable $A\subs X_f$ there is  
${\beta}\in {\omega}_2$ such that  $|A\cap \HH({\beta})|={\omega}_1$.
\end{lemma}

\proof
Assume that 
$p\force$ ``$\dot{A}=\{\dot{\alpha}_{\xi}:{\xi}<\oo\}\in\br \oot;\oo;$''.
For each ${\xi}<\oo$ pick $p_{\xi}\le p$ and ${\alpha}_{\xi}\in\oot$
such that $p_{\xi}\force\dot {\alpha}_{\xi}=\hat{\alpha}_{\xi}$.
Since $P_f$ satisfies c.c.c we can assume that the ${\alpha}_{\xi}$
are pairwise different. Let $\sup\{{\alpha}_{\xi}:{\xi}<\oo\}<{\beta}<\oot$.
Now for each ${\xi}<\oo$ define the condition $q_{\xi}\le p_{\xi}$
by the stipulations $a^{q_{\xi}}=a^{p_{\xi}}\cup\{{\beta}\}$, 
$h^{q_{\xi}}({\beta})=a^{q_{\xi}}$ and
 $i^{q_{\xi}}({\nu},{\beta})=\empt$ for ${\nu}\in a^{p_{\xi}}$.
Then $q_{\xi}\in P_f$ and $q_{\xi}\force\dot{\alpha}_{\xi}\in \HH({\beta})$.
But $P_f$ satisfies c.c.c, so there is $q\le p$ such that
$q\force$ ``$|\{{\xi}\in\oo:q_{\xi}\in\gcal\}|=\oo$''.
Thus $q\force$ ``$|A\cap\HH({\beta})|=\oo$.''
\eproof

Since every $\HH({\beta})$ is compact, lemma \ref{lm:oo-compact}
above clearly implies that $X_f$ is $\oo$-compact, i.e. every subset
$S\subs X_f$ of size $\oo$ has a complete accumulation point.

Now we start to work on (iii): in $V^{P_f}$ 
the closure of any countable subset $Y$ of 
$X_f$ is either compact or it contains a final segment of $\oot$. 
If $Y$ is also in the ground model, then actually the second alternative
occurs and this follows  easily from the next lemma.

\begin{lemma}
If $p\in P_f$, ${\beta}\in \ap$, $b\subs \ap\cap {\beta}$, 
${\alpha}\in {\beta}\setm \ap$, then there is a condition
$q\le p$ such that ${\alpha}\in \uuu q;{\beta};b;$.
\end{lemma}

\proof
Define the condition $q\le p$ by the following stipulations:
$\aq=\ap\cup \{{\alpha}\}$, $\hq({\alpha})=\{{\alpha}\}$,
$$
\hq({\nu})= \left\{ 
\begin{array}{ll}
\hp({\nu})\cup\{{\alpha}\}
&\mbox{if ${\beta}\in \hp({\nu})$}\\
\hp({\nu})
&\mbox{if ${\beta}\notin \hp({\nu})$}
\end{array}
\right.
$$
for ${\nu}\in\ap$, and let $\iq\supset\ip$
and $\iq({\alpha},{\nu})=\empt$ for ${\nu}\in\ap$.

To show $q\in P_f$  we 
need to check only \ref{df:poset}(iv). Assume that 
${\alpha}\in \hq({\nu})*\hq({\mu})$. Then  by the construction of $q$
we have 
${\beta}\in \hp({\nu})*\hp({\mu})$. Thus there is 
${\xi}\in \ip({\nu},{\mu})$ with ${\beta}\in \hp({\xi})$. But then
${\alpha}\in\hq({\xi})$, so by  $\iq({\nu},{\mu})=\ip({\nu},{\mu})$
we have ${\alpha}\in\hh q; \iq({\nu},{\mu});$. In view of
$\hp({\nu})*\hp({\mu})\subs \hh p;\ip({\nu},{\mu});$ we are done. 
Thus $q\in P_f$, $q\le p$ and clearly ${\alpha}\in  \uuu q;{\beta};b;$,
 so we are done. 
\eproof

This lemma  yields the following corollary.
\begin{corollary}\label{cl:old-set}
If $Z\in \br \oot;{\omega};\cap V$ and ${\beta}\in \oot\setm\sup Z$ 
then ${\beta}\in\overline{Z}$.
\end{corollary}

\proof
Let $\UU {\beta};b;$ be a neighbourhood of ${\beta}$, 
$b\in\br {\beta};<{\omega};$. Since 
$p\force $``$\UU {\beta};b;\supset\uuu p;{\beta};b;$'' for each 
$p\in P_f$ and the set
$$
D_{{\beta},b,Z}=\{q\in P_f: \uuu q;{\beta};b;\cap Z\ne\empt\}
$$
is dense in $P_f$ by the previous lemma, it follows that $\UU {\beta};b;$
intersects $Z$. Consequently ${\beta}\in\overline{Z}$.
\eproof

The space $X_f$ is right separated, i.e. scattered,  so we can consider its
Cantor-Bendixon hierarchy. According to corollary \ref{cl:old-set}
for each ${\alpha}<\oot$ the set
 $A_{\alpha}=[{\omega}{\alpha},{\omega}{\alpha}+{\omega})$ is a dense
set of isolated points  in
 $X_f\rest (\oot\setm {\omega}{\alpha})$. 
Thus the ${\alpha}^{\rm th}$  Cantor-Bendixon level of 
$X_f$ is just $A_{\alpha}$. Therefore $X_f$ is a thin very tall,
locally compact scattered space in the sense of \cite{Ro}. 
Let us emphasize that CH was not needed to get this
result, hence we have also given an alternative proof of the main result
of \cite{BS}.

Now we continue to work on proving  property \ref{tm:main}(iii) of $X_f$.

Given $p\in P_f$ and $b\subs \ap$ let $p\rest b=\<b,h,\ip\rest\br b;2; \>$
where $h$ is the function with $\dom(h)=b$ and $h({\xi})=\hp({\xi})\cap b$
for ${\xi}\in b$. Let us remark that $p\rest b$ in not necessarily in
 $P_f$. In fact,  $p\rest b\in P_f$ if and only if $\ip({\xi},{\eta})\subs b$
for each ${\xi}\ne {\eta}\in b$. Especially, if $b$ is an initial segment
of $\ap$, then $p\rest b\in P_f$. The order $\le$ of $P_f$ can be extended
in a natural way to the restrictions of conditions: if $p$ and $q$ are in
$P_f$, $b\subs \ap$, $c\subs \aq$, define $p\rest b\le q\rest c$ iff
$b\supset c$, $\hp({\xi})\cap c=\hq({\xi})\cap c$ for each ${\xi}\in c$,
and $\ip\rest \br c;2;=\iq\rest \br c;2;$. Clearly if $p\rest b\in P_f$ and 
$q\rest c\in P_f$ then the two definitions of $\le$ coincide.

\begin{definition}\label{df:blow}
Let $p,p'\in P_f$ with $a^p=a^{p'}$. 
We write  {\em $p\prec p'$} if for each ${\alpha}\in \ap$ and
$b\subs \ap\cap {\alpha}$ we have 
$\uuu p;{\alpha};b;\subset \uuu p';{\alpha};b;$.
\end{definition}

The following technical result will play a crucial role in the proof
of \ref{tm:main}(iii). Part (c) in it will enable us to ``insert''
certain things in $\HH({\gamma}_0)$ in a non-trivial way.
But there is a price we have to pay for this: this is the point
where the coherency of the $\HH({\alpha})$ has to be abandoned.
Part (d) will be needed in section \ref{sc:FU}.

\begin{lemma}\label{lm:amalg}
Assume that $s=\ahii s;\in P_f$, $\as=S\cup E\cup F$, $Q\subs S$,
$S \leo E \leo F$, $E=\{{\gamma}_i:i<k\}$, 
${\gamma}_0<{\gamma}_1<\dots <{\gamma}_{k-1}$,
$F=\{{\gamma}_{i,0},{\gamma}_{i,1}:i<k\}$, moreover  
\begin{enumerate}\rlabel
\item $\forall i<k$ $\hs({\gamma}_{i,0})\cap \hs({\gamma}_{i,1})=
\hhs{Q\cup E}$,
\item  $\forall i<k$ $\forall {\xi}\in S$
$f({\xi},{\gamma}_i)=f({\xi},{\gamma}_{i,0})=f({\xi},{\gamma}_{i,1})$.
\end{enumerate}
Then there is a condition $r=\ahii r;$ with  $\ar=S\cup E$  such that 
\begin{enumerate}\alabel
\item \label{a} $r\le s\rest S$,
\item \label{b} $r\le s\rest (Q\cup E)$, 
\item \label{c} $S\setm \hhs{Q\cup E}\subs\hr({\gamma}_0)$,
\item \label{d} $s\rest (S\cup E)\prec r$.
\end{enumerate}
\end{lemma}

\proof
Let $\ar=S\cup E$ and write $C=S\setm \hhs{Q\cup E}$. 
For ${\xi}\in\ar$ we set
\begin{displaymath}
\hr({\xi})= \left\{ 
\begin{array}{ll}
\hs({\xi})\cup C&
\mbox{if ${\xi}={\gamma}_i$ and ${\gamma}_0\in \hs({\gamma}_i)$,}\\
\hs({\xi})&\mbox{otherwise.}
\end{array}
\right.
\end{displaymath}

For ${\xi}\ne{\eta}\in\ar$ we let 
\begin{displaymath}
\ir({\xi},{\eta})= \left\{ 
\begin{array}{ll}
\is({\xi},{\eta})&\mbox{if ${\xi},{\eta}\in Q\cup E$ or ${\xi},{\eta}\in S$,}\\
f({\xi},{\eta})\cap \ar&\mbox{otherwise.}
\end{array}
\right.
\end{displaymath}
Finally let $r=\ahii r;$.

We claim that $r$ satisfies the requirements of the lemma.
\ref{a}, \ref{b} and \ref{c} are clear from the definition
of $r$, once we establish that $r\in P_f$.
To see that it suffices to
 check only \ref{df:poset}.(iv) because the other requirements are
clear from the construction of $r$. So let ${\xi}<{\eta}\in \ar$.
We have to show $\hr({\xi})*\hr({\eta})\subs \hhr {\ir({\xi},{\gamma}_i)}$.

If ${\xi},{\eta}\in S$, then 
$\hr({\xi})*\hr({\eta})\subs \hhr {\ir({\xi},{\eta})}$ 
holds because $r\rest S=s\rest S\in P_f$.

So we can assume that ${\eta}={\gamma}_i$ for some $i<k$.
\newcases
\begin{case}
${\xi}\in S$.
\end{case}
\begin{subcase}
${\xi}\notin \hr({\gamma}_i)$, hence $\hr({\xi})*\hr({\gamma}_i)=
\hr({\xi})\cap \hr({\gamma}_i)$.
\end{subcase}

In this case we also have $\hs({\xi})*\hs({\gamma}_i)=\hs({\xi})\cap \hs({\gamma}_i$
and so
\begin{equation}
\label{eq:*}
\hs({\xi})\cap\hs({\gamma}_i)\subs 
\hhs {\is({\xi},{\gamma}_i)}\subs \hhr {\ir({\xi},{\gamma}_i)}
\end{equation}
 for
$\ir({\xi},{\gamma}_i)\supset\is({\xi},{\gamma}_i,)$.
If ${\gamma}_0\notin \hs({\gamma}_i)$ then
$\hr({\gamma}_i)=\hs({\gamma}_i)$ and since $\hr({\xi})=\hs({\xi})$
we have $\hr({\xi})*\hr({\gamma}_i)=\hs({\xi})\cap \hs({\gamma}_i)\subs
\hhr {\ir({\xi},{\gamma}_i)}$ by (\ref{eq:*}). 
Assume now that ${\gamma}_0\in \hs({\gamma}_i)$. 
Thus $\hr({\gamma}_i) =\hs({\gamma}_i)\cup C$, and so ${\xi}\notin C$,
that is ${\xi}\in \hhs {Q\cup E}$. 
Then $\hr({\xi})*\hr({\gamma}_i)=\hr({\xi})\cap \hr({\gamma}_i)= 
(\hs({\xi})\cap \hs({\gamma}_i))\cup (\hs({\xi})\cap C)$. 
By (\ref{eq:*})
above it is enough to show that 
$\hs({\xi})\cap C\subs \hhr {\ir({\xi},{\gamma}_i)}$.
Since $\hs({\xi})\cap C=\empt$ for ${\xi}\in Q$ we can assume that 
${\xi}\notin Q$.
By (i) $\hhs {Q\cup E}=\hs({\gamma}_{i,0})\cap\hs({\gamma}_{i,1})$ and so
\begin{equation}
\label{eq:**}
\hs({\xi})\cap C=\hs({\xi})\setm \hhs{Q\cup E}=
%\hs({\xi})\setm (\hs({\gamma}_{i,0})\cap\hs({\gamma}_{i,1}))=
(\hs({\xi})\setm \hs({\gamma}_{i,0}))\cup 
(\hs({\xi})\setm \hs({\gamma}_{i,1})).
\end{equation}
Since ${\xi}\in \hhs{Q\cup E}\subs \hs({\gamma}_{i,j})$ for $j=0,1$
we have 
\begin{equation}\label{eq:x}
\hs({\xi})\setm \hs({\gamma}_{i,j})=
\hs({\xi})* \hs({\gamma}_{i,j})\subs \hhs {\is({\xi},{\gamma}_{i,j})}.
\end{equation}
By (ii), $\is({\xi},{\gamma}_{i,j})\subs
f({\xi},{\gamma}_{i,j})\cap \as=f({\xi},{\gamma}_i)\cap \as$.
Since ${\xi}\notin Q$ it follows that 
$f({\xi},{\gamma}_i)\cap \as= {\ir({\xi},{\gamma}_i)}$.
Thus from (\ref{eq:x}) we obtain
\begin{equation}\label{eq:xx}
\hs({\xi})\setm \hs({\gamma}_{i,j})\subs \hhs {\ir({\xi},{\gamma}_i)}
=\hhr {\ir({\xi},{\gamma}_i)}.
\end{equation}
Putting (\ref{eq:**}) and (\ref{eq:xx}) together we get 
$\hs({\xi})\cap C\subs \hhr {\ir({\xi},{\gamma}_i)}$ which was to be proved.

\begin{subcase}
${\xi}\in \hr({\gamma}_i)$,  hence $\hr({\xi})*\hr({\gamma}_i)=
\hr({\xi})\setm \hr({\gamma}_i)$.
\end{subcase}

If ${\xi}\in\hs({\gamma}_i)$ then since $\hs({\xi})=\hr({\xi})$ we have
$\hr({\xi})*\hr({\gamma}_i)=\hr({\xi})\setm \hr({\gamma}_i)\subs 
\hs({\xi})\setm \hs({\gamma}_i)\subs \hhs {\is({\xi},{\gamma}_i)}\subs
\hhr {\ir({\xi},{\gamma}_i)}$ and we are done.

So we can assume that ${\xi}\notin\hs({\gamma}_i)$ and so 
  $\hr({\gamma}_i)\ne \hs({\gamma}_i)$. By the construction of 
$r$,  we have
${\gamma}_0\in \hs({\gamma}_i)$,  
$\hr({\gamma}_i)=\hs({\gamma}_i)\cup C$ and so  ${\xi}\in C$, i.e 
${\xi}\notin \hhs{Q\cup E}$. By (i) we can assume that 
${\xi}\notin \hs({\gamma}_{i,0})$. 
So $s\in P_f$ implies
\begin{equation}\label{eq:zzz}
\hs({\xi})\cap \hs({\gamma}_{i,0})=
\hs({\xi})* \hs({\gamma}_{i,0})\subs
\hhs {\is({\xi},{\gamma}_{i,0})}
\end{equation}
We have
\begin{equation}\label{eq:z}
\hr({\xi})*\hr({\gamma}_i)=\hs({\xi})\setm (\hs({\gamma}_i)\cup C)\subs
\hs({\xi})\setm C 
\end{equation}
and applying (i) again 
\begin{equation}\label{eq:zz}
\hs({\xi})\setm C =
\hs({\xi})\cap \hhs{Q\cup E}\subs
\hs({\xi})\cap \hs({\gamma}_{i,0}).
\end{equation}
By (ii), $\is({\xi},{\gamma}_{i,0})\subs
f({\xi},{\gamma}_{i,0})\cap \as=f({\xi},{\gamma}_i)\cap \as$.
Since ${\xi}\notin Q\subs \hhs{Q\cup E}$ it follows that 
$f({\xi},{\gamma}_i)\cap \as= {\ir({\xi},{\gamma}_i)}$
and so (\ref{eq:zzz})--(\ref{eq:zz}) together yield
\begin{equation}
\hr({\xi})*\hr({\gamma}_i)\subs\hhr {\ir({\xi},{\gamma}_i)}
\end{equation}
which was to be proved.

\begin{case}
${\xi}={\gamma}_j$ for some $j<i$.
\end{case}

Since $\is({\gamma}_j,{\gamma}_i)=\ir({\gamma}_j,{\gamma}_i)$
we have
\begin{equation}\label{eq:xxx}
\hs({\gamma}_j)*\hs({\gamma}_i)\subs \hhr{\ir({\gamma}_i,{\gamma}_j)}.
\end{equation}
 It is easy to check that 
\begin{equation}\label{eq:xxxx}
\hr({\gamma}_j)*\hr({\gamma}_i)= \left\{ 
\begin{array}{ll}
\hs({\gamma}_j)*\hs({\gamma}_i)
&\mbox{if ${\gamma}_0\notin \hs({\gamma}_j)*\hs({\gamma}_i)$,}\\
\hs({\gamma}_j)*\hs({\gamma}_i)\cup C
&\mbox{if ${\gamma}_0\in \hs({\gamma}_j)*\hs({\gamma}_i)$.}
\end{array}
\right.
\end{equation}
So we are done if ${\gamma}_0\notin \hs({\gamma}_j)*\hs({\gamma}_i)$.
Assume ${\gamma}_0\in \hs({\gamma}_j)*\hs({\gamma}_i)$.
Then there is 
${\gamma}_l\in \is({\gamma}_j,{\gamma}_i)$ with 
${\gamma}_0\in\hs({\gamma}_l)$. Thus,
by the construction of $r$ we have
\begin{equation}\label{eq:xxxxx}
C\subs \hr({\gamma}_l)\subs \hhr{\ir({\gamma}_j,{\gamma}_i)}
\end{equation}
But (\ref{eq:xxxx}) and (\ref{eq:xxxxx}) together imply what we wanted.

Thus we proved $r\in P_f$.

Clearly $r$ satisfies \ref{lm:amalg}.(a)--(c). To check 
\ref{lm:amalg}.(d) write $s'=s\rest (S\cup E)$ and 
let ${\alpha}\in S\cup E$ and $b\subs (S\cup E)\cap {\alpha}$.
We need to show that $\uuu s';{\alpha};b;\subs \uuu r;{\alpha};b;$.
Since $S\cup E$ is an initial segment of $\as$, we have
$\uuu s';{\alpha};b;=\uuu s;{\alpha};b;$.
If ${\alpha}\in S$, then also $\uuu s';{\alpha};b;= \uuu r;{\alpha};b;$,
so we can assume ${\alpha}\in E$.

Let ${\xi}\in \uuu s';{\alpha};b;=\uuu s;{\alpha};b;$.
Then ${\xi}\in \hs({\alpha})\subs \hhs {E}$, and hence ${\xi}\notin C$.
But $\hhr {b}\setm \hhs{b}\subs C$, more precisely, it is empty or just $C$.
Since ${\xi}\in \uuu s;{\alpha};b;$, it follows that ${\xi}\notin \hhs {b}$
and so ${\xi}\notin \hhr{b}$ because ${\xi}\notin C$. Thus 
${\xi}\in \uuu r;{\alpha};b;$. Hence $r$ satisfies (d).

The lemma is proved.
\eproof

\begin{lemma}\label{lm:compact-or-tail}
In $V^{P_f}$, if $Y\subs \oot$ is countable, then either $\overline{Y}$ is
compact or $|\oot\setm\overline{Y}|\le\oo$.
\end{lemma}

\proof
Assume that $1_{P_f}\force$``$\ydot=\{\yydot_n:n\in{\omega}\}\subs \oot$''.
For each $n\in {\omega}$ fix a maximal antichain 
$C_n\subs P_f$ such that for each $p\in C_n$  there is
${\alpha}\in \ap$ with $p\force$``$\yydot_n=\hat{{\alpha}}$''.
Let $A=\bigcup\{\ap:p\in\bigcup\limits_{n<{\omega}}C_n\}$.
Since every $C_n$ is countable by c.c.c we have $|A|={\omega}$.

Assume also that $1_{P_f}\force$ 
``$\overline{\ydot}$ is not compact'', that is, 
$Y$ can not be covered by finitely many $H({\delta})$ in $V^{P_f}$.

Let
$$
I=\{{\delta}<\oot:\exists p\in P_f\ 
p\force\text{``${\delta}\notin \overline{\ydot}$''}\}.
$$
%$$
%I=\{{\delta}<\oot:\exists p_{\delta}\in P\ \exists D_{\delta}\in\fin {\delta};
%p_{\delta}\force \text{``$Y\cap(\UU {\delta};D_{\delta};)=\empt$''}\}.
%$$
Clearly 
$1_{P_f}\force\ \oot\setm I\subs\overline{\ydot}$.
Since $\oot\setm I$ is in the ground model,  
 by corollary \ref{cl:old-set}  it is enough to show that 
$\oot\setm I$ is infinite. Actually  we will prove much more:

\newClaim
\begin{Claim}
$I$ is not stationary in $\oot$.
\end{Claim}
%\Proof{the claim}
Assume on the contrary that $I$ is stationary.
Let us fix, for each ${\delta}\in I$, a condition $p_{\delta}\in P_f$ 
and a finite set $D_{\delta}\in\fin {\delta};$ such that 
$p_{\delta}\force$ ``$\ydot\cap\UU {\delta};D_{\delta};=\empt$''.
For each ${\delta}\in I$ let $Q_{\delta}=a^{p_{\delta}}\cap {\delta}$ and 
$E_{\delta}=a^{p_{\delta}}\setm {\delta}$.
We can assume that $D_{\delta}\subs Q_{\delta}$ and $\sup A <{\delta}$
for each ${\delta}\in I$.

Let $B_{\delta}=\clf(A\cup Q_{\delta},E_{\delta})$ for ${\delta}\in I$
(see \ref{lm:closure}). For each ${\delta}\in I$ the set  
$B_{\delta}$ is  countable  with 
$\sup(B_{\delta})=\sup(A\cup Q_{\delta})$, so  
 we can apply  Fodor's pressing down lemma and 
 CH  to get a stationary  set  $J\subs I$ and a countable set 
$B\subs \oot$ such that 
$B_{\delta}=B$ for each ${\delta}\in J$.

By thinning out $J$ and with a further use of  CH we can assume that
for a fixed $k\in {\omega}$ we have
\begin{enumerate}\arablabel
\item $E^{\delta}=\{{\gamma}^{\delta}_i:i<k\}$ for ${\delta}\in J$, 
${\gamma}^{\delta}_0<{\gamma}^{\delta}_1<\dots<{\gamma}^{\delta}_{k-1}$,
\item $f({\xi},{\gamma}^{\delta}_i)=f({\xi},{\gamma}^{{\delta}'}_i)$
for each ${\xi}\in B$, ${\delta},{\delta}'\in J$ and  $i<k$.
\end{enumerate}

Let ${\delta}=\min J$, $D=D_{\delta}$, $E=E_{\delta}$, $p=p_{\delta}$,
$Q=Q_{\delta}$.
By lemma \ref{lm:big_int} there are ordinals ${\delta}_j\in J$ with
${\delta}<{\delta}_0<{\delta}_1<\dots<{\delta}_{2k-1}$  such that
\begin{equation}\tag{$\star$}
B\cup E\subs \bigcap\{f({\xi},{\eta}):
{\xi}\in E_{{\delta}_i},{\eta}\in E_{{\delta}_j}, i<j<2k\}.
\end{equation}

For $i<k$ and $j<2$ let ${\gamma}_i={\gamma}^{\delta}_i$ 
and ${\gamma}_{i,j}={\gamma}^{{\delta}_{2i+j}}_i$. 
Let $F=\{{\gamma}_{i,j}:i<k,j<2\}$.

We know that $a^p=Q\cup E$. Define 
the condition $q\in P_f$ by the following
stipulations: \begin{enumerate}\rlabel
\item  $\aq=\ap\cup F$  and $q\le p$,
\item $\hq({\gamma}_{i,j})=\{{\gamma}_{i,j}\}\cup \ap$
for $\<i,j\>\in k\times 2$
\item $\iq({\gamma}_{i_0,j_0}, {\gamma}_{i_1,j_1})=\ap$
for $\<i_0,j_0\>\ne \<i_1,j_1\>\in k\times 2$,
\item $\iq({\xi},{\gamma}_{i,j})=\empt$ for ${\xi}\in\ap$ and
 $\<i,j\>\in k\times 2$.
\end{enumerate}
Since $\ap\subs B\cup E$, ($\star$) implies that $q\in P_f$.

Since $1_{P_f}\force$ ``$\HHH {Q\cup E}\not\supset\ydot$'',  
there is a condition $t\le q$, 
a natural number $n$ and an ordinal ${\alpha}$ such that 
$t\force$ ``${\alpha}=\yydot_n$'' but
${\alpha}\in\at\setm \hht {Q\cup E}$. Since $C_n$ is a maximal antichain 
we can assume that $t\le v$ for some $v\in C_n$.

Let $s=t\rest(B\cup E\cup F)$. Then $s\in P_f$ 
because for each pair ${\xi}<{\eta}\in\as$ if ${\xi}\in B$ then 
$\itt({\xi},{\eta})\subs f({\xi},{\eta})\subs B$ and so 
$\itt({\xi},{\eta})\subs \as$, and if ${\xi},{\eta}\in E\cup F$ then
$\itt({\xi},{\eta})=\iq({\xi},{\eta})\subs Q\cup E\subs \as$.
Moreover $s\le v$ because $a^{v}\subs B$. Thus 
$s\force$ ``${\alpha}=\yydot_n$'' and ${\alpha}\notin\hhs{Q\cup E}$. Let
$S=\as\cap B$. 

Since $\is({\gamma}_{i,0}, {\gamma}_{i,1})=
\iq({\gamma}_{i,0}, {\gamma}_{i,1})=Q\cup E$,  and 
${\gamma}_{i,j}\notin \hs({\gamma}_{i,1-j})$ we have
\begin{equation}\label{eq:sups}
\hs({\gamma}_{i,0})\cap\hs({\gamma}_{i,1})=
\hs({\gamma}_{i,0})*\hs({\gamma}_{i,1})\subs \hhs{Q\cup E}.
\end{equation}
Moreover, if ${\xi}\in Q\cup E$ and $j<2$  
then ${\xi}\in\hq ({\gamma}_{i,j})\subs
\hs({\gamma}_{i,j})$ and  $\is({\xi},{\gamma}_{i,j})=\empt$, 
consequently 
\begin{equation}\label{eq:subs}
\hs({\xi})\setm \hs({\gamma}_{i,j})=\hs({\xi})*\hs({\gamma}_{i,j})\subs
\hhs{\is({\xi},{\gamma}_{i,j})}=\empt.
\end{equation}
Putting (\ref{eq:sups}) and (\ref{eq:subs})  together it follows that
$\hs({\gamma}_{i,0})\cap\hs({\gamma}_{i,1})=\hhs{Q\cup E}$.

Thus we can apply \ref{lm:amalg} to get a condition $r$ such that 
$r\le s\rest (Q\cup E)=p$, $r\le s\rest S \le v$ and 
${\alpha}\in S\setm \hhs{Q\cup E}\subs \hs({\gamma}_0)=\hs({\delta})$.
Since $D\subs Q$, we have ${\alpha}\in \hs({\delta})\setm \hhs D$.

Thus
$$
r\force {\alpha}\in \ydot\cap (H({\delta})\setm \HHH D).
$$
On the other hand $r\force \ydot\cap (H({\delta})\setm \HHH D)=\empt$
because $r\le p$. With this contradiction the claim  is proved and this
completes the proof of the lemma.
%\Eproof
%The lemma is proved.
\eproof

Clearly, lemma \ref{lm:compact-or-tail} implies that $X_f$ is 
countably compact.

\begin{corollary}\label{lm:oszt}
If $F\subs X$ is closed (or open), 
then either $|F|\le\oo$ or $|X\setm F|\le\oo$.
\end{corollary}

\proof
If $|F|=\oot$ then $F$ is not compact, so by lemma
 \ref{lm:noncompact_free_sequence} $F$ contains a free sequence $Y$
with non-compact closure.
But $F$ is  initially $\oo$-compact and countably tight, 
so  $Y$ is countable. Consequently,  
we have $|\oot\setm\overline{Y}|\le\oo$ by lemma \ref{lm:compact-or-tail} 
and so $|X\setm F|\le\oo$.
\eproof

\begin{corollary}\label{lm:zx}
$X_f$ is normal and $\z(X_f)=\operatorname{hd}(X_f)\le \oo$.
\end{corollary}

\proof
To show that $X_f$ is normal let
$F_0$ and $F_1$ be disjoint closed subsets of $X_f$.  Since at least 
one of them is compact by lemma \ref{lm:compact-or-tail} they
 can be separated by open subsets of $X_f$ because
 $X_f$ is $T_3$.

Concerning the hereditarily density of $X_f$
it follows easily from corollary \ref{lm:oszt} that 
$X_f$ does not contain a discrete subspace of size
$\oot$. But $X_f$ is right separated, 
so all the left separated subspaces of $X_f$  are of size 
$\le \oo$, that is , $\z(X)\le\oo$.
\eproof

Thus theorem \ref{tm:main} is proved.
\Eproof

We know that the space $X_f$ is not automatically hereditarily
separable, so the following question   of Arhangel'ski\u\i,  
\cite[problem 5]{Arh}, remains unanswered: 
Is it true in ZFC that every hereditarily separable, 
initially $\oo$-compact space is compact?

As we have seen our space $X_f$ is normal. However, we don't know
whether $X_f$ is or can be made   hereditarily normal, i.e. $T_5$. 
This raises the following problem.
\begin{problem}
Is it provable in ZFC that every $T_5$, countably tight, initially
$\oo$-compact space is compact? 
\end{problem}

\section{Making $X_f$ Frechet-Uryson}\label{sc:FU}

In  \cite[problem 12]{Arh}  Arhangel'ski\u\i asks if it is provable in ZFC that a 
normal, first countable initially $\oo$-compact space is 
necessarily  compact. We could not completely answer this 
question, but in this section we show that the Frechet-Uryson property
(which is sort of half-way between countable tightness and first countability)
in not enough to get compactness.

To achieve that we want to find a further extension of 
the model $V^{P_f}$ in which 
 $X_f$ becomes Frechet-Uryson but   its other properties are preserved, 
for example, $X_f$ remains initially $\oo$-compact and normal.
Since $X_f$ is countably tight and ${\chi}(X_f)\le\oo$
it is a  natural idea to make 
$X_f$  Frechet-Uryson by constructing a 
generic extension of $V^{P_f}$ in which  $X_f$ remains countably tight
and ${\rm p}>\oo$, i.e. 
${\rm MA}_{\oo}(\text{${\sigma}$-centered})$ holds 
(see \cite[theorem 8]{Sho}).

The standard c.c.c poset $P$ which forces $p>\oo$
is obtained by a suitable finite support iteration   
of length $2^{\oo}$. During this iteration in the  ${\alpha}^{\rm th}$ step
we choose a non-principal filter $\fcal\subs \pcal({\omega})$ generated
by at most $\oo$ elements and we add a new subset $A$ of ${\omega}$ to the 
${\alpha}^{\rm th }$ intermediate model so that $A$ is almost contained in
every element of $\fcal$, i.e. $A\setm F$ is finite for each $F\in\fcal$. 
 It is well-known and easy to see that  $P$   has property $K$. 
Thus, by   theorem \ref{tm:preservation_omega_tight} 
below,
$X_f$ remains  countably tight  in $V^{P_f*R}$ and so indeed 
$X_f$ becomes Frechet-Uryson in that model. 
Moreover,  theorem \ref{tm:preservation_oo-compact} implies that
the $\oo$-compactness of $X_f$ is also preserved.  Unfortunately,
we could not prove that forcing with $P$ 
preserves the  countable compactness of $X_f$.   

So, instead of aiming at  $p>\oo$  we will consider only
those filters during the iteration 
which  are needed  in proving the Frechet-Uryson property of
$X_f$. As we will see, we can handle these filters in such a way that 
our iterated forcing $R$ preserves not only the countable 
compactness of $X_f$ but also property (iii) from theorem
\ref{tm:main}: in $V^{P_f*R}$ the closure of any countable subset  of $X_f$ is either compact or   contains a final segment of $\oot$. 
Of course, this will insure the preservation of the normality
of $X_f$ as well.

We start  with the two easy theorems, promised above,  about the
preservation of countable tightness and 
$\oo$-compactness of $X_f$ under certain c.c.c
forcings.

\begin{theorem}
\label{tm:preservation_omega_tight}
If the topological space $X$ is right separated, compact, countably tight
and the poset $R$ has property $K$ 
then forcing with $R$ preserves the  countable tightness of $X$.
\end{theorem}

\proof
First we recall that $X$ remains compact (and clearly right separated) 
in any extension of the ground model by \cite[lemma 7]{JW}.
Since $\operatorname{F}(X)=\tght(X)$
for compact spaces, 
assume  indirectly that    
$1_R\force$ ``{\em $\{\zzdot_{\xi}:{\xi}<\oo\}\subs X$ is a free sequence}''. 
For every ${\xi}<\oo$ we have that
$1_R\force$``{\em $\overline{\{\zzdot_{\zeta}:{\zeta}<{\xi}\}}$ and
$\overline{\{\zzdot_{\zeta}:{\xi}\le {\zeta}<\oo\}}$ 
are disjoint compact sets}'' and $X$ is $T_3$, 
so we can fix  a condition $p_{\xi}\in P$,
 open sets $U_{\xi}$ and 
$V_{\xi}$ from the ground model and a point $z_{\xi}\in X$ such that 
$\overline{U_{\xi}}\cap \overline{V_{\xi}}=\empt$ and
$$
p_{\xi}\force \text{``
 $\{\zzdot_{\zeta}:{\zeta}<{\xi}\}\subs U_{\xi}$,  
$\{\zzdot_{\zeta}:{\xi}\le {\zeta}<\oo\}\subs V_{\xi}$
and $\zzdot_{\xi}=z_{\xi}$.''}
$$ 
 
Since $R$ has property K, there is an uncountable set $I\subs \oo$
such that the conditions $\{p_{\xi}:{\xi}\in I\}$ are pairwise compatible.

We claim that the sequence $\{z_{\xi}:{\xi}\in I\}$ is an uncountable free 
sequence in the ground model which 
contradicts $\operatorname{F}(X)=\tght(X)={\omega}$.
Indeed let ${\xi}\in I$. If ${\zeta}\in I\cap {\xi}$, then $p_{\xi}$
and $p_{\zeta}$ has a common extension $q$ in $P$ and we have
$$
q\force \text{``$\zzdot_{\zeta}=z_{\zeta}$  and
$\{\zzdot_{\eta}:{\eta}<{\xi}\}\subs \udot_{\xi}$.''}
$$
Hence $z_{\zeta}\in U_{\xi}$. Similarly for ${\zeta}\in I\setm {\xi}$
we have $z_{\zeta}\in V_{\xi}$. Therefore
$\overline{U_{\xi}}$ and $\overline{V_{\xi}}$ separate 
$\overline{\{z_{\zeta}:{\zeta}\in I\cap {\xi}\}}$
and $\overline{\{z_{\zeta}:{\zeta}\in I\setm {\xi}\}}$ which implies that
$\{z_{\xi}:{\xi}\in I\}$ is really free.
\eproof

\begin{theorem}\label{tm:preservation_oo-compact} 
Forcing with a c.c.c poset $R$ over $V^{P_f}$ preserves property
 \ref{tm:main}.$($ii$)$ of the space $X_f$, i.e.
for each uncountable $A\subs X_f$ there is 
${\beta}\in \oot$ such that $A\cap H({\beta})$ is uncountable.
\end{theorem}

\proof
We work in $V^{P_f}$.
Assume that $r\force_{R}$ 
``$\adot=\{\dot{{\alpha}}_{\xi}:{\xi}<\oo\}\in\br X_f;\oo;$''.
For each ${\xi}<\oo$ pick a condition $r_{\xi}\le r$ from $R$ which
decides the value of $\dot{{\alpha}}_{\xi}$,
$r_{\xi}\force_{R}$ ``$\dot{{\alpha}}_{\xi}={\alpha}_{\xi}$''. 
Since $R$ satisfies 
c.c.c, $\{{\alpha}_{\xi}:{\xi}\in\oo\}$ is uncountable, hence as 
$X_f$ has property (ii) in $V^{P_f}$, for some ${\beta}<\oot$
the set  $I=\HH({\beta})\cap \{{\alpha}_{\xi}:{\xi}<\oo\}$ is 
also uncountable.
Since $R$ satisfies c.c.c there is a condition $q\le r$ in $R$
such that $q\force_{R}$ ``$|\{{\xi}\in I:r_{\xi}\in \gcal\}|=\oo$'',
where $\gcal$ is the $R$-generic filter over $V^{P_f}$.
Thus $q\force_{R}$ ``$|\adot\cap \HH({\beta})|=\oo$''
which was to be proved.
\eproof

Of course, theorem \ref{tm:preservation_oo-compact} implies that 
forcing with any c.c.c poset $R$  preserves the $\oo$-compactness of $X_f$.
It is  much harder to find a property of  a poset $R$ which
guarantees that forcing with $R$ over $V^{P_f}$
preserves the countable compactness of $X_f$.
We will proceed in the following way.
In definition \ref{df:nice} we formulate when a poset $\rdot$ is called  
{\em nice (over $P_f$)}, and then in theorem 
\ref{tm:spec-preservation} we show that forcing with a nice poset 
preserves not only the countable compactness of $X_f$, but also 
property \ref{tm:main}(iii):  the closure of any countable subset  of $X_f$ 
is either compact or  contains a final segment of $\oot$. Finally,
in definitions \ref{df:FU-poset} and \ref{df:FU-iteration} we describe 
a class of finite support iterated forcings, which 
by theorem \ref{tm:FU-iteration} are nice and have property K , 
and then in theorem \ref{tm:FU} we show that forcing with a 
suitable member of this class makes $X_f$ Frechet-Uryson.

\begin{definition}\label{df:nice}
Let $\rdot$ be a name for a poset in $V^{P_f}$. We say that $\rdot$
is {\em nice (over $P_f$)} if there is a dense
subset $\dcal$ of the iteration $P_f*\rdot$ with the following property:
%\begin{enumerate}\nolabel%\item 
If $\<p_0,r_0\>, \<p_1,r_1\>\in \dcal$,
$p,p'\in P_f$  are such that $p\le p_0,p_1$, $p'\le p_0,p_1$, and 
$$
p\force  \text{``$r_0$ and $r_1$ are compatible in $\rdot$'',}
$$
moreover we have  either $p\le p'$ or $p\prec p'$ (see definition 
\ref{df:blow})  then we also have 
$$p'\force \text{``$r_0$ and $r_1$ are compatible in $\rdot$''.}
$$
%\end{enumerate}
\end{definition}

\begin{theorem}\label{tm:spec-preservation}
If CH holds in $V$, $f$ is a $\Delta$-function and
$\rdot$ is a $P_f$-name for a  c.c.c poset
which is nice over $P_f$, then \ref{tm:main}.$($iii$)$ is preserved
by forcing with $\rdot$, i.e.
$$
V^{P_f*\rdot}\models\forall Y\in\br X_f;{\omega};
(\ \overline{Y} 
\text{ is compact or }|X_f\setm \overline{Y}|\le \oo\ ).
$$
\end{theorem}

\proof
Let $\dcal\subs P_f*\rdot$ witness that $\rdot$ is nice.

Assume that $1_{P*\rdot}\force$``$\ydot=\{\yydot_n:n\in{\omega}\}\subs \oot$''.
For each $n\in {\omega}$ fix a maximal antichain 
$C_n\subs \dcal$ such that for each $\<p,r\>\in C_n$  there is
${\alpha}\in \ap$ with $\<p,r\>\force$``$\yydot_n=\hat{{\alpha}}$''.
Let $A=\bigcup\{\ap:\<p,r\>\in\bigcup\limits_{n<{\omega}}C_n\}$.
Since every $C_n$ is countable by c.c.c we have $|A|={\omega}$.

Assume  that $1_{P_f*\rdot}\force$ 
``$\overline{\ydot}$ is not compact'', that is, 
$\ydot$ can not be covered by finitely many $H({\delta})$ in $V^{P_f*\rdot}$.

Let
$$
I=\{{\gamma}<\oot:\exists \<p,r\>\in P_f*\rdot\ 
\<p,r\>\force\text{``${\gamma}\notin \overline{\ydot}$''}\}.
$$

Since $\oot\setm I$ is in the ground model and 
$1_{P_f*\rdot}\force\ \oot\setm I\subs\overline{\ydot}$,  
  it is enough to show that $\oot\setm I$ is infinite.
Indeed, in this case $\overline{\oot\setm I}$  
contains a final segment of $\oot$ by corollary \ref{cl:old-set}. 
(Note that the closure of a ground model set does not change 
under any further forcing.)
Thus the next claim completes the proof of this theorem.

\newClaim
\begin{Claim}
$I$ is not stationary.
\end{Claim}

Assume on the contrary that $I$ is stationary.
For each  ${\delta}\in I$  fix  a condition 
$\<p_{\delta},r_{\delta}\>\in P_f*\rdot$ and a finite set 
$D_{\delta}\in\fin {\delta};$ such that 
$\<p_{\delta},r_{\delta}\>\force$ ``$Y\cap\UU {\delta};D_{\delta};=\empt$''.
For each ${\delta}\in I$ write $Q_{\delta}=a^{p_{\delta}}\cap {\delta}$ and 
$E_{\delta}=a^{p_{\delta}}\setm {\delta}$.
We can assume that $D_{\delta}\subs Q_{\delta}$, ${\delta}\in E_{\delta}$
and $\sup A<{\delta}$ for each ${\delta}\in I$.

Let $B_{\delta}=\clf(A\cup Q_{\delta},E_{\delta})$ for ${\delta}\in I$
(see \ref{lm:closure}). Since 
$B_{\delta}$ is a countable set with 
$\sup(B_{\delta})=\sup(A\cup Q_{\delta})$ and so $\sup(B_{\delta})<{\delta}$,
we can apply  Fodor's pressing down lemma and 
CH to get a stationary  set  $J\subs I$ and a countable set 
$B\subs \oot$ such that $B_{\delta}=B$ for each ${\delta}\in J$.

By thinning out $J$ and with a another use of CH  we can assume that 
for some
fixed $k\in {\omega}$ we have
\begin{enumerate}\arablabel
\item $E^{\delta}=\{{\gamma}^{\delta}_i:i<k\}$ for ${\delta}\in J$, 
${\gamma}^{\delta}_0<{\gamma}^{\delta}_1<\dots<{\gamma}^{\delta}_{k-1}$,
\item $f({\xi},{\gamma}^{\delta}_i)=f({\xi},{\gamma}^{{\delta}'}_i)$
for each ${\xi}\in B$, ${\delta},{\delta}'\in J$ and  $i<k$.
\end{enumerate}

Let ${\delta}=\min J$, $D=D_{\delta}$, $E=E_{\delta}$, $p=p_{\delta}$,
$r=r_{\delta}$,  $Q=Q_{\delta}$.
By lemma \ref{lm:big_int} there are $2k$ ordinals ${\delta}_j\in J$ with
${\delta}<{\delta}_0<{\delta}_1<\dots<{\delta}_{2k-1}$ from $J$ such that
\begin{equation}\tag{$\star$}
B\cup E\subs \bigcap\{f({\xi},{\eta}):
{\xi}\in E_{{\delta}_i},{\eta}\in E_{{\delta}_j}, i<j<2k\}.
\end{equation}

For $i<k$ and $j<2$ let ${\gamma}_i={\gamma}^{\delta}_i$ 
and ${\gamma}_{i,j}={\gamma}^{{\delta}_{2i+j}}_i$. 
Let $F=\{{\gamma}_{i,j}:i<k,j<2\}$.

We know that $a^p=Q\cup E$. Define the condition $q\in P_f$ by the following
stipulations: \begin{enumerate}\rlabel
\item $\aq=\ap\cup F$ and $q\le p$, 
\item $\hq({\gamma}_{i,j})=\{{\gamma}_{i,j}\}\cup \ap$
for $\<i,j\>\in k\times 2$
\item $\iq({\gamma}_{i_0,j_0}, {\gamma}_{i_1,j_1})=\ap$
for $\<i_0,j_0\>\ne \<i_1,j_1\>\in k\times 2$,
\item $\iq({\xi},{\gamma}_{i,j})=\empt$ for ${\xi}\in\ap$ and
 $\<i,j\>\in k\times 2$.
\end{enumerate}
Since $\ap\subs B\cup E$, 
($\star$) implies that $q\in P_f$ and so $\<q,r\>\in P*\rdot$. 

Since $1_{P_f*\rdot}\force$ ``{\em $\HHH {Q\cup E}$ does not cover $Y$}'', 
 there is a condition $\<t,u\>\le \<q,r\>$ in $P_f*\rdot$, 
a natural number $n$ and an ordinal ${\alpha}$ such that 
$\<t,u\>\force$ ``${\alpha}=\yydot_n$'' but
${\alpha}\in \at\setm \hht {Q\cup E}$. Since $C_n$ is a maximal antichain 
we can assume that $\<t,u\>\le \<c,d\>$ for some $\<c,d\>\in C_n$.

Let $s=t\rest(B\cup E\cup F)$, $w=t\rest (S\cup E)$ and 
$S=a^s\cap B=a^w\cap B$. 
Then $s$ and $w$ are in $P_f$ 
because for each pair ${\xi}<{\eta}\in\as$ if ${\xi}\in B$ then 
$\itt({\xi},{\eta})\subs f({\xi},{\eta})\subs B$ and so 
$\itt({\xi},{\eta})\subs S$, and if ${\xi},{\eta}\in E\cup F$ then
$\itt({\xi},{\eta})=\iq({\xi},{\eta})\subs Q\cup E\subs \as\cap a^w$.
Moreover $s\le w\le p,c$ because $a^{c}\subs a^t\cap B=S$ and 
$\ap=Q\cup E\subs \at\cap\subs B\cup E$. 

Since $\<t,u\>\le \<p,r\>$ and $\<t,u\>\le \<c,d\>$, so
\begin{equation}\tag{*}
t\force \text{`` $r$ and $d$ are compatible in $\rdot$''}.
\end{equation}
Since $t\le w\le p,c$ and $\rdot$ is nice for $P_f$ we have  
\begin{equation}\tag{**}
w\force \text{`` $r$ and $d$ are compatible in $\rdot$''}.
\end{equation}
We know $s\le c$ and so
$s\force $``${\alpha}=\yydot_n$'' and ${\alpha}\in\as\setm\hhs{Q\cup E}$.
 
Since $\is({\gamma}_{i,0}, {\gamma}_{i,1})=
\iq({\gamma}_{i,0}, {\gamma}_{i,1})=Q\cup E$,  and 
${\gamma}_{i,j}\notin \hs({\gamma}_{i,1-j})$ we have
\begin{equation}\label{eq:sups2}
\hs({\gamma}_{i,0})\cap\hs({\gamma}_{i,1})=
\hs({\gamma}_{i,0})*\hs({\gamma}_{i,1})\subs \hhs{Q\cup E}.
\end{equation}
If ${\xi}\in Q\cup E$, then ${\xi}\in\hq ({\gamma}_{i,j})\subs
\hs({\gamma}_{i,j})$ and  $\is({\xi},{\gamma}_{i,j})=\empt$ 
and so
\begin{equation}\label{eq:subs2}
\hs({\xi})\setm \hs({\gamma}_{i,j})=\hs({\xi})*\hs({\gamma}_{i,j})\subs
\hhs{\is({\xi},{\gamma}_{i,j})}=\empt.
\end{equation}
Putting (\ref{eq:sups2}) and (\ref{eq:subs2})  together it follows that
$\hs({\gamma}_{i,0})\cap\hs({\gamma}_{i,1})=\hhs{Q\cup E}$.
%Since $\iq({\gamma}_{i,0}, {\gamma}_{i,1})=\ap$,
% $\iq({\xi},{\gamma}_{i,j})=\empt$ for ${\xi}\in\ap$ and 
%$\hq({\gamma}_{i,0})\cap\hq({\gamma}_{i,1})=Q\cup E$
%it follows that 
%$\hs({\gamma}_{i,0})\cap\hs({\gamma}_{i,1})=\hs(\ap)=\hhs{Q\cup E}$.
Thus we can apply \ref{lm:amalg} to get a condition $p^*\in P_f$ such that 
$p^*\le s\rest (Q\cup E)=p$, $p^*\le s\rest S \le c$,  $w\prec p^*$ and
${\alpha}\in S\setm \hhs {Q\cup E}\subs h^{p^*}({\gamma}_0)=h^{p^*}({\delta})$.
But $D\subs Q$, so ${\alpha}\in h^{p^*}({\delta})\setm \hh p^*;D;$.
Hence $p^*\force $``${\alpha}\in \UU{\delta};D;$''.

Since $p^*\le p, c$ and $w\prec p^*$ and $\rdot$ is nice, it follows 
from (**) that
\begin{equation}\tag{***}
p^*\force \text{`` $r$ and $d$ are compatible in $\rdot$''}.
\end{equation}
Thus  $\<p,r\>$ and $\<c,d\>$ have    a common extension
$\<p^*,r^*\>$ in $P*\rdot$.

But then
$$
\<p^*,r^*\>\force \yydot_n={\alpha}\in \ydot\cap \UU {\delta}; D;
$$
because $\<p^*,r^*\>\le \<c,d\>$.
On the other hand $\<p^*,r^*\>\force$ ``$\ydot\cap \UU {\delta};D;=\empt$''
because $\<p^*,r^*\>\le \<p,r\>$. Contradiction,
the claim is proved, which completes the proof of the theorem.
\eproof

\begin{definition}\label{df:FU-poset}
Assume that  $A\in\br X_f;{\omega};$ and ${\alpha}\in\overline{A}$.
We define the poset  $Q(A,{\alpha})$ as follows. Its underlying set is
$\fin A;\times \fin {\alpha};$ . If $\<s,C\>$ and $\<s'C'\>$ are 
conditions, let $\<s,C\>\le \<s'C'\>$ if and only if
$s\supset s'$, $C\supset C'$ and $s\setm s'\subs \UU {\alpha};C';$.
For $q\in Q(A,{\alpha})$ write $q=\<s^q,C^q\>$ and $\supp(q)=s^q\cup C^q$.
\end{definition}

It is well-known and easy to see that if 
$\gcal$ is a $Q(A,{\alpha})$-generic filter, then
$S=\bigcup\{s^q:q\in\gcal\}$ is a sequence from $A$ which converges to 
${\alpha}$, i.e. every open neighbourhood  of ${\alpha}$ contains
all but finitely many points of $S$.

Clearly $Q(A,{\alpha})$ is ${\sigma}$-centered and well-met.
In fact, if $p_0=\<s_0,C_0\>$ and $p_1=\<s_1,C_1\>$ are compatible,
then $p_0\land p_1=\<s_0\cup s_1, C_0\cup C_1\>$.

\begin{definition}\label{df:FU-iteration}
A finite support  iterated forcing $\<R_{\xi}:{\xi}\le {\kappa}\>$ 
over $V^{P_f}$ is called an {\em FU-iteration} if 
for each ${\xi}<{\kappa}$ we have
$$
1_{P_f*R_{\xi}}\force 
\text{$R_{{\xi}+1}=R_{\xi}*Q(\adot,\dot{\alpha})$ for some
$\adot\in\br X_f;{\omega};$ and $\dot{\alpha}\in \overline{\adot}$}.
$$
\end{definition}

Since every $Q(A,{\alpha})$ is ${\sigma}$-centered, it is clear that
any FU-iteration is c.c.c, in fact it even has property K. 
The really important, and much less trivial, property of them
is given in our next result.

\begin{theorem}\label{tm:FU-iteration}
Any FU-iteration is nice over $P_f$.
\end{theorem}

\proof
Assume that $\<R_{\xi}:{\xi}\le {\kappa}\>$ is an FU-iteration,
$R_{{\xi}+1}=R_{\xi}*Q(\adot_{\xi},\dot{\alpha}_{\xi})$.
Write $Q^*=\fin \oot;\times \fin \oot;$. Clearly 
$1_{P_f*R_{\xi}}\force$ ``$Q(\adot_{\xi},\dot{\alpha}_{\xi})\subs Q^*$''.

We consider the elements of $P_f*\rdot_{\kappa}$ 
as  pairs $\<p,\rrdot\>$, where
$p\in P_f$ and $p\force$ ``{\em $\rrdot$ is a finite function, 
$\dom(\rrdot)\subs {\kappa}$ and 
$\rrdot\rest {\xi}\force_
{\rdot_{\xi}} \rrdot({\xi})\in Q(\adot_{\xi},\dot{\alpha}_{\xi})$
for each ${\xi}\in\dom (\rrdot)$}''.
\begin{definition}\label{df:determined}
A condition $\<p,r\>\in P_f*R_{\kappa}$ is called {\em determined}
if $r\in\fn({\kappa},Q^*,{\omega})$, (i.e $r$ is a finite function, not only
a $P_f$-name of a finite function)  moreover 
for each ${\xi}\in \dom (r)$
we have $\supp(r({\xi}))\subs \ap$ and
there is an ordinal ${\alpha}\in\ap$ such that
$\<p,r\rest {\xi}\>\force $ ``$\dot{\alpha}_{\xi}=\hat{\alpha}$''.
%Write $\supp_p(r)=\{{\alpha}(p,r,{\xi}):{\xi}\in\dom(r)\}\cup\bigcup
%\{\supp(r({\xi}):{\xi}\in \dom(r)\}$.
\end{definition}
Clearly the family $\dcal$ of the determined conditions is dense in
$P_f*R_{\kappa}$. We claim that $\dcal $  
witnesses that $R_{\kappa}$ is nice.
Indeed let $\<p_0,r_0\>, \<p_1,r_1\>\in \dcal$,
$p,p'\in P_f$  be such that $p\le p_0,p_1$, $p'\le p_0,p_1$, 
\begin{equation}
\tag{+}p\force  \text{``$r_0$ and $r_1$ are compatible in $\rdot$''}
\end{equation}
and either $p\le p'$ or $p\prec p'$.

Write $r_i({\xi})=\<s_i({\xi}),C_i({\xi})\>$ for ${\xi}\in\dom(r_i)$.
Let $D=\dom(r_0)\cup\dom(r_1)$. 
For $i<2$ and ${\xi}\in D$ let
\begin{displaymath}
s^*_i({\xi})= \left\{ 
\begin{array}{ll}
s_i({\xi})&\mbox{if ${\xi}\in\dom(r_i)$,}\\
\empt&\mbox{otherwise,}
\end{array}
\right.
\end{displaymath}
and
\begin{displaymath}
C^*_i({\xi})= \left\{ 
\begin{array}{ll}
C_i({\xi})&\mbox{if ${\xi}\in\dom(r_i)$,}\\
\empt&\mbox{otherwise.}
\end{array}
\right.
\end{displaymath}

Consider the condition  $r^*\in\fn({\kappa},Q^*,{\omega})$ defined by the stipulations $\dom(r^*)=D$ and
$r^*({\xi})=\<s^*_0({\xi})\cup s^*_1({\xi}),C^*_0({\xi})\cup C^*_1({\xi})\>$
for ${\xi}\in D$.
We show that $p'\force$ ``{\em $r^*$ is a common extension of
$r_0$ and $r_1$ in $R_{\kappa}$.}''
We prove a bit more:  we show, by a finite induction, 
that for each  ${\xi}\in D$
$$
p'\force \text{``$r^*\rest {\xi}+1\in R_{{\xi}+1}$ and $r^*\rest {\xi}+1\le
r_0\rest {\xi}+1, r_1\rest {\xi}+1$.''}
$$

The non-trivial step is when ${\xi}\in \dom(r_0)\cap \dom(r_1)$. 

Since  $\<p',r^*\rest {\xi}\>\le \<p',r_i\rest {\xi}\>$ 
by the induction hypothesis, it follows that
 $\<p',r^*\rest {\xi}\>\force$ 
``$ s_i({\xi})\in \br \adot_{\xi};<{\omega};$'' and so
$\<p',r^*\rest {\xi}\>\force$ 
``$r^*({\xi})\in Q(\adot_{\xi},\dot{{\alpha}_{\xi}}) $'', that is,
$p'\force$ ``$r^*\rest {\xi}+1\in R_{{\xi}+1}$''.

By (+) for some $P_f$-name $\rrdot$ and   $P_f* R_{\xi}$-name $\qqdot$ 
we have $\<p,\rrdot\>\in P_f* R_{\xi}$ and
$\<p,\rrdot\>\force$ 
``{\em $\qqdot\in Q(\adot_{\xi},\dot{\alpha}_{\xi})$ 
is a common extension of $r_0({\xi})$ and $r_1({\xi})$}''.
Hence $\<p,\rrdot\>\force $ ``$s^{\qqdot}\setm s_0({\xi})\subs \UU{\alpha};C_0({\xi});$ 
and $s^{\qqdot}\supset s_1({\xi})$''.
Therefore $\<p,\rrdot\>\force$ 
$``s_1({\xi})\setm s_0({\xi})\subs \UU{\alpha};C_0({\xi});$''.
By $\supp(r_i({\xi}))\subs \ap$, this can only happen if 
\begin{equation}
\tag{\dag}
s_1({\xi})\setm s_0({\xi})\subs \uuu p; {\alpha}; C_0({\xi});.
\end{equation}
But both $p\le p'$ and $p\prec p'$ together with (\dag) imply
\begin{equation}
\tag{\ddag}
s_1({\xi})\setm s_0({\xi})\subs \uuu p'; {\alpha}; C_0({\xi});.
\end{equation}
Thus $p'\force$ 
``$s_1({\xi})\setm s_0({\xi})\subs \UU {\alpha};C_0({\xi});$'',
i.e. $\<p', r^*\rest {\xi}\>$ ``$\force r^*({\xi})\le r_0({\xi})$''.
Similarly, it can be proved that 
 $\<p', r^*\rest {\xi}\>\force $ ``$r^*({\xi})\le r_1({\xi})$''.
Thus we have carried out the induction step and  the theorem is proved.
\eproof

Now we are ready to prove the main result of this section.
\begin{theorem}\label{tm:FU}
Assume that CH holds in the ground model $V$, 
there is a  $\Delta$-function $f$
 and ${\lambda}$ is a cardinal
such that $\oo<{\lambda}={\lambda}^{\omega}$. Then in $V^{P_f}$ there
is an FU-iteration $\rdot_{\lambda}$ of length ${\lambda}$ 
such that in $V^{P_f*\rdot_{\lambda}}$
the space $X_f$ is Frechet-Uryson and satisfies 
\ref{tm:main}$($i$)$--$($iii$)$, 
moreover 
$V^{P_f*\rdot_{\lambda}}\models$``$2^{\kappa}=({\lambda}^{\kappa})^V$ 
for each cardinal ${\kappa}\ge {\omega}$''.
\end{theorem}

\proof
Since $|P_f|=\oot$ and $P_f$ satisfies c.c.c we have
$({\lambda}^{\omega})^{V^{P_f}}\le
 ((|P_f|{\lambda})^{\omega})^{V}=({\lambda}^{\omega})^{V}={\lambda}$.
Therefore,  using a suitable book-keeping procedure in $V^{P_f}$
(see \cite[Ch VIII. 6.3]{Ku} for this technique)
 we can construct an  FU-iteration  
$\<R_{\xi}:{\xi}\le {\lambda}\>$, 
$R_{{\xi}+1}=R_{\xi}*Q(\adot_{\xi},\dot{\alpha}_{\xi})$,
 having the following property: for every pair $\<A,{\alpha}\>$ if
 $A$  is a countable subset of $X_f$
in $V^{P_f*R_{\lambda}}$ and ${\alpha}\in\overline{A}$ then for some
${\xi}<{\lambda}$ we have 
$V^{P_f*R_{\xi}}\models$ ``$Q(A_{\xi},{\alpha}_{\xi})=Q(A,{\alpha})$''.
Thus 
\begin{enumerate}
\item[(*)] 
$V^{P_f*R_{\lambda}}\models$ ``{\em if ${\alpha}\in\oot$ 
is in the closure of a countable set $A\subs X_f$ 
then there is a sequence  $\{s_n:n\in {\omega}\}\subs A$ which converges
to ${\alpha}$}''. 
\end{enumerate} 
Since $R_{\lambda}$ has property K  in $V^{P_f}$, by theorem \ref{tm:preservation_omega_tight}
the space $X_f$ remains countably tight in $V^{P_f*R_{\lambda}}$.
Putting together this observation with (*) it follows that  
$V^{P_f*R_{\lambda}}\models $ ``$X_f$ is Frechet-Uryson''.

An FU-iteration is nice by theorem \ref{tm:FU-iteration},
and so $X_f$ has property \ref{tm:main}(iii) in $V^{P_f*R_{\lambda}}$
by theorem \ref{tm:spec-preservation}.
Since $X_f$ is right separated it remains locally compact in any 
extension of $V^{P_f}$.
Since $R_{\lambda}$ satisfies c.c.c the space $X_f$
has property \ref{tm:main}(ii)  in $V^{P_f*R_{\lambda}}$ by theorem
\ref{tm:preservation_oo-compact}. All this implies that $X_f$
remains initially $\oo$-compact and normal.

Finally we investigate the cardinal exponents in $V^{P_f*\rdot_{\lambda}}$.
Since $|P_f|=\oot$ and $1_{P_f}\force$``$|R_{\lambda}|={\lambda}$'',  
the iteration $P_f*R_{\lambda}$ contains a dense subset $\dcal$ of
cardinality $\le{\lambda}$. Since $P_f*R_{\lambda}$ satisfies c.c.c
 it follows that for each ${\kappa}\ge {\omega}$ we have
 $(2^{\kappa})^{V^{P_f*R_{\lambda}}}\le ((|\dcal|^{\omega})^{\kappa})^{V}=
({\lambda}^{\kappa})^V$.

On the other hand every successor step of
an FU-iteration introduces a new subset of  a countable set,  and so
$(2^{\omega})^{V^{P_f*R_{\lambda}}}\ge {\lambda}$. Consequently,
${V^{P_f*R_{\lambda}}}\models$``$2^{\kappa}=
(2^{{\omega})^{\kappa}}\ge {\lambda}^{\kappa}\ge ({\lambda}^{\kappa})^V
$'', which proves what we wanted.
\eproof

Theorem \ref{tm:FU} answers a question raised  by Arhangel'ski\u\i, 
\cite[problem 3]{Arh}, in the negative: CH can not be
 weakened to $2^{\omega}<2^{\oo}$ in the theorem of van Douwen and Dow.
In fact we proved much more:
the existence of a   Frechet-Uryson, initially ${\omega}_1$-compact and 
non-compact space  is consistent with practically any cardinal arithmetic 
that violates CH. More precisely, if we have a ZFC model $V$ in which
CH holds and ${\lambda}={\lambda}^{\omega}\le 2^{\oo}$, then we can find 
a cardinal preserving generic extension $W$ of $V$ which contains a
Frechet-Uryson and normal counterexample   to the van Douwen--Dow question, 
$(2^{\omega})^W={\lambda}$, moreover $(2^{\kappa})^W=(2^{\kappa})^V$
for each ${\kappa}\ge \oo$.
We can obtain $W$ as follows. First
we force the ${\sigma}$-complete poset $P$ of Shelah  (see \cite{BS}) which 
introduces a $\Delta$-function $f$  in $V^P$. 
Since $P$ is ${\sigma}$-complete and $|P|=\oot$, 
forcing with $P$ does not change $2^{\kappa}$ for any 
${\kappa}\ge {\omega}$.
Now forcing with
$P_f$ over $V^P$ introduces the counterexample $X_f$ to the
van Douwen--Dow question. Since $|P_f|=\oot$, the cardinal exponents
are the same in $V^P$ and in $V^{P_f}$ for uncountable cardinals and
$(2^{\omega})^{V^{P*P_f}}=\oot$. Finally we can apply theorem
\ref{tm:FU} to get the desired final model $W=V^{P*P_f*R_{\lambda}}$.
Let us remark that  we have $(2^{\kappa})^W=(2^{\kappa})^V$ 
for ${\kappa}\ge \oo$ because
$(2^{\kappa})^V=({\lambda}^{\kappa})^V$ by ${\lambda}\le 2^{\oo}$.

Let us remark that for any cardinal ${\kappa}$ 
the poset $\fn({\kappa},2,{\omega})$,
i.e. the forcing notion that adds ${\kappa}$ many  Cohen reals, 
is clearly nice over $P_f$, as is witnessed by the dense set
of the determined conditions.
Thus, by theorems \ref{tm:preservation_omega_tight},
\ref{tm:preservation_oo-compact} and
\ref{tm:spec-preservation} , adding Cohen reals will preserve
properties \ref{tm:main}(i)--(iii) of $X_f$, especially $X_f$ remains 
 countably tight and initially $\oo$-compact.  
It is  worthwhile to mention that, in contrast with this, 
Alan Dow proved in \cite{Dow} that if  CH  holds in the ground model
 $V$ then adding Cohen reals can not introduce
 a  countably tight, initially $\oo$-compact and  non-compact $T_3$ 
space.

Let us finish by formulating the following higher cardinal version of the van Douwen--Dow problem:
\begin{problem}
Is it provable in ZFC
that an initially $\oot$-compact $T_3$ space of countable tightness is 
compact ?
\end{problem}

The main problem is trying to use out approach 
that worked for $\oot$ (instead of ${\omega}_3$) here is that 
no $\Delta$-function may exist for ${\omega}_3$!
This problem has come up already in the efforts trying to lift the result
of \cite{BS} from $\oot$ to ${\omega}_3$.

\end{document}